\documentclass{article}
\usepackage{latexsym}
\usepackage[varumlaut]{yfonts}
\usepackage{epsfig}
\usepackage{amsmath}
\usepackage{amssymb}
\usepackage{amsthm}
\newtheorem{theorem}{Theorem}[section]
\newtheorem{corollary}[theorem]{Corollary}

\newtheorem{lemma}[theorem]{Lemma}

\theoremstyle{remark}
\newtheorem{remark}[theorem]{Remark}

\theoremstyle{definition}
\newtheorem{definition}[theorem]{Definition}

\DeclareMathOperator\diag{diag}

\DeclareMathOperator\pf{pf}

\begin{document}

\title{Centro-affine hypersurface immersions with parallel cubic form}

\author{Roland Hildebrand \thanks{%
LJK, Universit\'e Grenoble 1 / CNRS, 51 rue des Math\'ematiques, BP53, 38041 Grenoble cedex, France
({\tt roland.hildebrand@imag.fr}).}}

\maketitle

\begin{abstract}
We consider non-degenerate centro-affine hypersurface immersions in $\mathbb R^n$ whose cubic form is parallel with respect to the Levi-Civita connection of the affine metric. There exists a bijective correspondence between homothetic families of proper affine hyperspheres with center in the origin and with parallel cubic form, and K\"ochers conic $\omega$-domains, which are the maximal connected sets consisting of invertible elements in a real semi-simple Jordan algebra. Every level surface of the $\omega$ function in an $\omega$-domain is an affine complete, Euclidean complete proper affine hypersphere with parallel cubic form and with center in the origin. On the other hand, every proper affine hypersphere with parallel cubic form and with center in the origin can be represented as such a level surface. We provide a complete classification of proper affine hyperspheres with parallel cubic form based on the classification of semi-simple real Jordan algebras. Centro-affine hypersurface immersions with parallel cubic form are related to the wider class of real unital Jordan algebras. Every such immersion can be extended to an affine complete one, whose conic hull is the connected component of the unit element in the set of invertible elements in a real unital Jordan algebra. Our approach can be used to study also other classes of hypersurfaces with parallel cubic form.
\end{abstract}



\section{Introduction}

The cubic form $C$ of an equiaffine hypersurface immersion is the covariant derivative of the affine metric $h$ with respect to the affine connection $\nabla$. Affine hypersurface immersions with parallel cubic form have been studied in various settings for more than 20 years, and their classification is an important problem in affine differential geometry.

One can consider immersions whose cubic form is parallel with respect to the affine connection, $\nabla C = 0$. Non-degenerate Blaschke immersions satisfying this condition are either quadrics or graph immersions whose graph function is a cubic polynomial \cite{Vrancken88}. Actually, in the latter case the immersion must be an improper affine hypersphere \cite{BNS90}, and the determinant of the Hessian of the graph function identically equals $\pm1$ \cite[p.47]{NomizuSasaki}. Non-degenerate Blaschke hypersurface immersions with $\nabla C = 0$ into $\mathbb R^k$, $k = 3,4,5,6$, have been classified in \cite{NomizuPinkall89},\cite{Vrancken88},\cite{Gigena02},\cite{Gigena03}, respectively. In \cite{Gigena11} an algorithm was presented to classify all such immersions for a given arbitrary dimension.

Another class of hypersurface immersions is obtained when the condition $\nabla K = 0$ is imposed, where $K = \nabla - \hat\nabla$ is the difference tensor between the affine connection and the Levi-Civita connection of the affine metric. Non-degenerate Blaschke immersions with this property have been studied in \cite{DillenVrancken98}. There it was established that, as for $\nabla C = 0$, they are either quadrics or improper affine hyperspheres. In the latter case the graph function is given by a polynomial, the affine metric is flat, the difference tensor is nilpotent, i.e., $K_X^m = 0$ for some $m > 1$ and all vector fields $X$, and $[K_X,K_Y] = 0$ for all vector fields $X,Y$.

Parallelism of the cubic form can also be defined with respect to the connection $\hat\nabla$. Since the affine metric is parallel with respect to $\hat\nabla$, the conditions $\hat\nabla C = 0$ and $\hat\nabla K = 0$ are equivalent. Much work concentrated on the case of Blaschke immersions. A Blaschke immersion satisfying $\hat\nabla C = 0$ must be an affine hypersphere \cite{BNS90}. In \cite{MagidNomizu89} all Blaschke immersions into $\mathbb R^3$ satisfying $\hat\nabla C = 0$ have been classified. In \cite{DillenVrancken91},\cite{HuLi11},\cite{HLLV11b} all such Blaschke immersions into $\mathbb R^4$ with definite, Lorentzian, and general affine metric, respectively, have been classified. In \cite{DVY94} all such Blaschke immersions into $\mathbb R^5$ with definite affine metric have been classified, and it has been shown that in arbitrary dimension, definiteness of the affine metric implies that the immersion is either a quadric or a locally homogeneous affine hypersphere. In \cite{HLSV09} all such immersions into $\mathbb R^k$, $k \leq 8$, with definite affine metric have been classified. In \cite{DillenVrancken94},\cite{HLV08} it has been observed that the Calabi product of affine hyperspheres with parallel cubic form or of such a hypersphere with a point are again affine hyperspheres with parallel cubic form, and hence one can speak of decomposable or irreducible such immersions. In a classification, one then only needs to consider the irreducible immersions. Finally, in \cite{HLV11},\cite{HLLV11a} a classification of all irreducible Blaschke hypersurface immersions with parallel cubic form whose affine metric is definite or Lorentzian, respectively, has been achieved.

A closer look at the classification in \cite{HLV11} reveals that the locally strongly convex hyperbolic affine hyperspheres with parallel cubic form are exactly those hyperspheres which are asymptotic to symmetric cones. Now the interiors of the symmetric cones are exactly the convex $\omega$-domains of K\"ocher \cite{Koecher99}, and it is not hard to verify that the hyperspheres in question are exactly the level surfaces of the $\omega$-function in these $\omega$-domains (we will give a brief account on $\omega$-domains in Subsection \ref{sec_inv}). One of the main results of this contribution is that this relation holds in general, i.e., independently of the convexity assumption. Namely, every non-degenerate proper affine hypersphere with center in the origin satisfying $\hat\nabla C = 0$ can be represented as a level surface of the $\omega$-function in some $\omega$-domain (Theorem \ref{th_AHS2}), and conversely, every such level surface is a non-degenerate proper affine hypersphere with center in the origin satisfying $\hat\nabla C = 0$ (Theorem \ref{S1}). Since the $\omega$-function is homogeneous \cite[p.35]{Koecher99}, it is clear that every ray of the $\omega$-domain intersects the affine hypersphere exactly once. We may then define a projection $\pi$ from the $\omega$-domain onto the hypersphere, taking every ray to its intersection point with the hypersphere.

The $\omega$-domains of K\"ocher are closely linked to real semi-simple Jordan algebras $J$. Namely, every $\omega$-domain can be represented as a connected component of the set of invertible elements in $J$, and every such connected component is an $\omega$-domain. The classification of proper affine hyperspheres with parallel cubic form then reduces to the classification of real semi-simple Jordan algebras. Much like in the case of semi-simple Lie algebras, any semi-simple Jordan algebra breaks down into a direct sum of a finite number of simple algebras, each of which is in turn a member of one of finitely many infinite series, or one of finitely many exceptional algebras. We show that the Calabi product of affine hyperspheres with parallel cubic form corresponds to the decomposition of semi-simple Jordan algebras into simple factors (Lemma \ref{lem_Calabi}). This allows to characterize the proper affine hyperspheres with parallel cubic form as Calabi products of irreducible such hyperspheres. The irreducible proper affine hyperspheres with parallel cubic can in turn be classified using the classification of simple real Jordan algebras (Theorem \ref{main_sphere}).

However, the main subject of this contribution are centro-affine hypersurface immersions satisfying $\hat\nabla C = 0$. In \cite{LiuWang97} the centro-affine hypersurface immersions into $\mathbb R^3$ satisfying $\hat\nabla \tilde C = 0$ have been classified, where $\tilde C$ is the traceless part of the cubic form. Since $\hat\nabla C = 0$ implies $\hat\nabla \tilde C = 0$, this classification includes also the immersions satisfying $\hat\nabla C = 0$. In \cite[Theorem 1.3]{LiWang91} the centro-affine hypersurface immersions satisfying $\hat\nabla C = 0$ with flat definite affine metric have been classified. It turns out that these are, up to centro-affine equivalence, the surfaces given by
\begin{equation} \label{LiWangFamily}
\prod_{k=1}^n x_k^{\alpha_k} = 1,\qquad \alpha_k > 0.
\end{equation}
If one relaxes the flatness condition to positive or negative semi-definiteness of the Ricci tensor, then the immersion may also be a proper affine hypersphere \cite[Prop. 7.2]{LiLiSimon04}. Note that a proper affine hypersphere with center in the origin and with parallel cubic form has also parallel cubic form as a centro-affine immersion. Thus the classifications in \cite{HLV11},\cite{HLLV11a} yield also examples of centro-affine immersions with parallel cubic form. It has been observed in \cite[p.342]{LiLiSimon04} that the family \eqref{LiWangFamily} of centro-affine hypersurfaces with parallel cubic form contains exactly one affine hypersphere. Moreover, it is not hard to see that all hypersurfaces in this family are asymptotic to the boundary of the same $\omega$-domain, namely the interior of the positive orthant, and that these hypersurfaces are coverings of the hypersphere under the projection $\pi$ of the $\omega$-domain onto the hypersphere.

We show that many of these properties hold for centro-affine hypersurface immersions with parallel cubic form in general. The role of the $\omega$-domains is played by the connected components $Y$ of the set of invertible elements in real unital Jordan algebras. This class of algebras is wider than that of the real semi-simple Jordan algebras. However, an analog $\varphi$ of the $\omega$-function can still be defined on $Y$. Every centro-affine hypersurface immersion satisfying $\hat\nabla C = 0$ is a homogeneous symmetric space, contained in some connected component $Y$, and is a covering of the level surfaces of $\varphi$ under an analog of the projection $\pi$, with the different sheets of the covering being homothetic images of each other (Theorem \ref{th_main3}). On the other hand, every unital real Jordan algebra satisfying an extra condition, namely the existence of a non-degenerate {\sl trace form} $\gamma$ \cite[p.24]{Schaefer}, defines centro-affine hypersurfaces with parallel cubic form (Theorem \ref{alg_to_imm}). The classification of centro-affine hypersurface immersions with parallel cubic form is hence reduced to the classification of real unital Jordan algebras with non-degenerate trace form. The decomposition of a unital Jordan algebra into a direct product of lower-dimensional algebras does not in general generate a decomposition of the corresponding centro-affine immersions into lower-dimensional ones. The generalization of the Calabi product proposed in \cite[Example 3.4]{LiLiSimon04} turns out to be not sufficient to capture all cases of centro-affine immersions with decomposable Jordan algebra.

The remainder of the paper is structured as follows. In the next section we introduce a convenient way to work with centro-affine hypersurface immersions, namely by describing them as integral manifolds of an involutive distribution. In Section \ref{sec_Jordan} we provide the necessary background on Jordan algebras. Section \ref{sec_correspondence} contains the main technical results of the paper, namely how exactly unital Jordan algebras are related to centro-affine hypersurface immersions with parallel cubic form. In Section \ref{sec_classification} we achieve the classification of those centro-affine hypersurface immersions with parallel cubic form which correspond to semi-simple Jordan algebras. This includes a complete classification of proper affine hyperspheres with parallel cubic form. In Sections \ref{sec_desc},\ref{sec_correspondence},\ref{sec_classification}, whenever we speak about {\it parallel cubic form}, we will mean the condition $\hat\nabla C = 0$ for a centro-affine hypersurface immersion. In Section \ref{sec_other} we apply our methods to other classes of hypersurface immersions with parallel cubic form, in particular, those mentioned at the beginning of this introduction.

\section{A description of centro-affine immersions} \label{sec_desc}

In this section we introduce a description of smooth non-degenerate centro-affine hypersurface immersions $f: M \to \mathbb R^n$ by logarithmically homogeneous functions and reformulate the condition $\hat\nabla C = 0$ in terms of these functions. Under more restrictive conditions, we will also use a description of centro-affine immersions by integral manifolds of an involutive distribution. The distribution will in turn be described as the kernel of a closed 1-form, with the logarithmically homogeneous functions being its local potentials. This will allow us to establish a link with Jordan algebras in Section \ref{sec_correspondence}.

\medskip

Let $f: M \to \mathbb R^n$ be a smooth centro-affine hypersurface immersion. Consider the direct product ${\cal D} = M \times \mathbb R_{++}$ of $M$ with the open half-line of positive real numbers. Let us define an immersion $\tilde f: {\cal D} \to \mathbb R^n$ by
\begin{equation} \label{tilde_f}
\tilde f: (\xi,\lambda) \mapsto \lambda f(\xi).
\end{equation}
Since $f$ is a centro-affine immersion, the immersion $\tilde f$ will be a local diffeomorphism. By means of the differential $\tilde f_*$ we can then identify the tangent space $T_y{\cal D}$ with the target space $\mathbb R^n$ for every $y \in {\cal D}$. In particular, we can define on ${\cal D}$ the position vector field $e$ (the reason for this notation will become apparent in Section \ref{sec_correspondence}), such that $e$ at the point $(\xi,\lambda) \in {\cal D}$ is given by $(\tilde f_*)^{-1}(\lambda f(\xi))$. Moreover, ${\cal D}$ inherits from $\mathbb R^n$ its flat affine connection $D$.

We shall now define a function $F: {\cal D} \to \mathbb R$. For $(\xi,\lambda) \in {\cal D}$ we set
\begin{equation} \label{F_def}
F(\xi,\lambda) = \log\lambda.
\end{equation}
Thus $F$ is {\sl logarithmically homogeneous} of degree $\nu = 1$, i.e., it satisfies the relation
\begin{equation} \label{log_homogeneity}
F(\xi,\alpha\lambda) = \nu\log\alpha + F(\xi,\lambda)
\end{equation}
with $\nu = 1$.

We shall adopt the Einstein summation convention over repeating indices. We will denote the derivatives of $F$ with respect to the flat affine connection $D$ by indices after a comma. Thus in an affine coordinate chart $y^{\alpha}$ on ${\cal D}$ we have $\frac{\partial F}{\partial y^{\alpha}} = F_{,\alpha}$, $\frac{\partial^2 F}{\partial y^{\alpha}\partial y^{\beta}} = F_{,\alpha\beta}$ etc.

We adopt the convention that the transversal vector field on the centro-affine immersion equals minus the position vector. Thus locally strongly convex immersions of hyperbolic type will have a negative definite affine metric $h$, while those of elliptic type will have a positive definite affine metric.

{\lemma \label{FderX} Let $f: M \to \mathbb R^n$ be a smooth centro-affine hypersurface immersion, let the function $F: {\cal D} = M \times \mathbb R_{++} \to \mathbb R$ be defined by \eqref{F_def}, and denote by $e$ the position vector field on ${\cal D}$. Then we have
\[ F_{,\alpha}e^{\alpha} = 1,\quad F_{,\alpha\beta}e^{\beta} = -F_{,\alpha},\quad F_{,\alpha\beta\gamma}e^{\gamma} = -2F_{,\alpha\beta},\quad F_{,\alpha\beta\gamma\delta}e^{\delta} = -3F_{,\alpha\beta\gamma},\quad F_{,\alpha\beta}e^{\alpha}e^{\beta} = -1.
\] }

\begin{proof}
The first relation is obtained by differentiating \eqref{log_homogeneity} with respect to $\alpha$ at $\alpha = 1$. The next three relations are obtained by differentiating repeatedly the first relation. The fifth relation follows from the first two.
\end{proof}

Note that
\begin{equation} \label{F10}
F_{,\alpha}u^{\alpha} = 0
\end{equation}
for every vector $u \in T_y{\cal D}$ which is tangent to the fiber $M \times \{\lambda\}$ at $y = (\xi,\lambda)$.

The following result, which is a generalization of \cite[Theorem 1, p.428]{Loftin02a}, then links the centro-affine pseudo-metric $h$ on $M$ with the Hessian metric $F''$ on ${\cal D}$.

{\theorem  \label{splitting} Let $f: M \to \mathbb R^n$ be a smooth centro-affine hypersurface immersion and let the function $F: {\cal D} = M \times \mathbb R_{++} \to \mathbb R$ be defined by \eqref{F_def}. Then the Hessian $F'' = D^2F$ on ${\cal D}$ is the orthogonal sum of the centro-affine form $h$ on $M$ and the form $-\lambda^{-2}d\lambda^2$ on $\mathbb R_{++}$.

In particular, if the immersion $f$ is non-degenerate, then the Hessian $F''$ is non-degenerate, and ${\cal D}$, seen as a pseudo-Riemannian manifold equipped with the Hessian pseudo-metric $F''$, equals the direct product $(M,h) \times (\mathbb R_{++},-\lambda^{-2}d\lambda^2)$ of pseudo-Riemannian manifolds, where $h$ is the centro-affine pseudo-metric on $M$. }

\begin{proof}
Let $\sigma(t)$ be a geodesic of the affine connection $\nabla$ on $M$. Consider the curve $\gamma(t) = (\sigma(t),1)$ on ${\cal D}$. Then the second derivative $\ddot\gamma$ is proportional to the transversal vector field $-e$, and the proportionality factor by definition equals the value of the affine fundamental form $h$ on the vector $\dot\gamma$. In an affine coordinate chart on ${\cal D}$ we have $\dot F = F_{,\alpha}\dot\gamma^{\alpha}$,
\[ \ddot F = F_{,\alpha\beta}\dot\gamma^{\alpha}\dot\gamma^{\beta} + F_{,\alpha}\ddot\gamma^{\alpha} = F_{,\alpha\beta}\dot\gamma^{\alpha}\dot\gamma^{\beta} - F_{,\alpha\beta}e^{\beta}\ddot\gamma^{\alpha} = F_{,\alpha\beta}\dot\gamma^{\alpha}\dot\gamma^{\beta} + F_{,\alpha\beta}e^{\beta}h(\dot\gamma,\dot\gamma)e^{\alpha} = F_{,\alpha\beta}\dot\gamma^{\alpha}\dot\gamma^{\beta} - h(\dot\gamma,\dot\gamma),
\]
where the second and fourth equality comes from the second and fifth relation in Lemma \ref{FderX}, respectively. Since $F \equiv 0$ on the curve $\gamma$, we obtain $h(\dot\gamma,\dot\gamma) = F_{,\alpha\beta}\dot\gamma^{\alpha}\dot\gamma^{\beta}$. Hence the restriction of the symmetric bilinear form $F''$ on the submanifold $M \times \{1\} \subset {\cal D}$ equals the affine fundamental form $h$ on $M$.

By \eqref{log_homogeneity} the tensor field $F''$ is invariant with respect to homotheties of ${\cal D}$, i.e., mappings of the form $(\xi,\lambda) \mapsto (\xi,\alpha\lambda)$ for fixed $\alpha > 0$. Hence the restriction of $F''$ to any level surface of $F$ equals $h$.

From the fifth relation in Lemma \ref{FderX} it follows that for every $\xi \in M$, the ray $\{\xi\} \times \mathbb R_{++} \subset {\cal D}$, equipped with the restriction of the symmetric bilinear form $F''$, is a Riemannian space with (negative definite) metric $-\lambda^{-2}d\lambda^2$. Finally, by \eqref{F10} and the second relation in Lemma \ref{FderX} the position vector field $e$ is orthogonal to the level surfaces of $F$ with respect to $F''$. This completes the proof.
\end{proof}

We consider only non-degenerate immersions $f$ and can hence use the pseudo-metric $F_{,\alpha\beta}$ to raise and lower indices of tensors on ${\cal D}$. Denote the elements of the inverse of the Hessian $F''$ by $F^{,\alpha\beta}$. Raising the index $\alpha$ in the second relation of Lemma \ref{FderX}, and contracting with $e^{\beta}$ and using the first relation of Lemma \ref{FderX} yields
\begin{equation} \label{FderX2}
F^{,\alpha\beta}F_{,\beta} = -e^{\alpha},\quad F^{,\alpha\beta}F_{,\alpha}F_{,\beta} = -1,
\end{equation}
respectively.

Let $\Pi: {\cal D} \to M$ be the projection onto $M$. For $W$ a covariant tensor field of order $k$ on $M$, the pullback $\Pi^*W$ is a covariant tensor field of the same order $k$ on ${\cal D}$. The pullbacks $\Pi^*W$ will allow us to work with tensors defined on the source manifold $M$ of the centro-affine immersion while using the very convenient affine coordinate system of the target space $\mathbb R^n$. The next result shows how different centro-affine invariants can be represented by expressions depending on the derivatives of $F$.

{\lemma \label{correspondence} Let $f: M \to \mathbb R^n$ be a smooth non-degenerate centro-affine hypersurface immersion and let the function $F: {\cal D} = M \times \mathbb R_{++} \to \mathbb R$ be defined by \eqref{F_def}. Consider the pullback of covariant tensors on $M$ to ${\cal D}$ by the projection $\Pi: {\cal D} \to M$. We have that
\begin{itemize}
\item the centro-affine pseudo-metric $h$ maps to $\Pi^*h_{\alpha\beta} = F_{,\alpha\beta} + F_{,\alpha}F_{,\beta}$.
\item the cubic form $C$ maps to $\Pi^*C_{\alpha\beta\gamma} = F_{,\alpha\beta\gamma} + 2F_{,\alpha\beta}F_{,\gamma} + 2F_{,\alpha\gamma}F_{,\beta} + 2F_{,\beta\gamma}F_{,\alpha} + 4F_{,\alpha}F_{,\beta}F_{,\gamma}$.
\item the Tchebycheff form $T = tr_hC$ maps to $\Pi^*T_{\alpha} = F_{,\alpha\beta\gamma}F^{,\beta\gamma} + 2nF_{,\alpha}$.
\item the covariant derivative $\hat\nabla C$ of the cubic form with respect to the Levi-Civita connection of $h$ maps to $\Pi^*\hat\nabla_{\delta}C_{\alpha\beta\gamma} = F_{,\alpha\beta\gamma\delta} - \frac12F^{,\rho\sigma}(F_{,\alpha\beta\rho}F_{,\gamma\delta\sigma} + F_{,\alpha\gamma\rho}F_{,\beta\delta\sigma} + F_{,\alpha\delta\rho}F_{,\beta\gamma\sigma})$.
\end{itemize} }

\begin{proof}
First note that the expressions on the right-hand sides in Lemma \ref{correspondence} are fully symmetric. With the relations in Lemma \ref{FderX} and \eqref{FderX2} it is easily verified that the contractions of these expressions with the position vector $e$ evaluate to zero. Therefore the values of these tensors on tangent vectors $u,v,\dots$ depend only on the projections $\Pi_*(u),\Pi_*(v),\dots$ of these tangent vectors.

Let $u,v,w$ be vector fields on ${\cal D}$ which are tangent to the fibers $M \times \{\lambda\}$.

We have $h(\Pi_*u,\Pi_*v) = F''(u,v) = F''(u,v) + F'(u)F'(v)$ by Theorem \ref{splitting} and \eqref{F10}, which proves the first item.

The cubic form $C$ on $M$ is the derivative of the centro-affine pseudo-metric $h$ with respect to the centro-affine connection $\nabla$. By definition of $\nabla$ we have $\nabla_{\Pi_*u}(\Pi_*v) = D_uv + \Pi^*h(u,v)e$. Hence we have
\[ C(\Pi_*u,\Pi_*v,\Pi_*w) = (D_wF'')(u,v) - F''(v,e)\Pi^*h(u,w) - F''(u,e)\Pi^*h(v,w) = F'''(u,v,w).
\]
Here the second relation comes from the fact that $e$ is orthogonal to $u,v$ in the Hessian pseudo-metric $F''$. The second item now follows from \eqref{F10}.

By Theorem \ref{splitting} the pullback of the Tchebycheff form $T$ is the trace of the pullback $\Pi^*C$ of the cubic form with respect to the Hessian pseudo-metric $F''$, $\Pi^*T = tr_{F''}\Pi^*C$. Inserting the expression of $\Pi^*C$ given by the second item of Lemma \ref{correspondence} and using the second relation in \eqref{FderX2}, we prove the third item.

The Christoffel symbols of the Levi-Civita connection $\hat D$ of the Hessian pseudo-metric $F''$ on ${\cal D}$ are given by
\begin{equation} \label{Christoffel}
\Gamma_{\alpha\beta}^{\gamma} = \frac12F^{,\gamma\delta}F_{,\alpha\beta\delta}.
\end{equation}
The covariant derivative of $F'$ then equals
\[ \hat D_\alpha F_{,\beta} = F_{,\alpha\beta} - \frac12F_{,\gamma}F^{,\gamma\delta}F_{,\alpha\beta\delta} = F_{,\alpha\beta} + \frac12e^{\delta}F_{,\alpha\beta\delta} = 0,
\]
where the second equality follows from the first relation in \eqref{FderX2}, and the third equality from the third relation in Lemma \ref{FderX}. Hence $F'$ is parallel with respect to $\hat D$. From the second item in Lemma \ref{correspondence} it then follows that the covariant derivative of $\Pi^*C$ with respect to $\hat D$ equals the covariant derivative of $F'''$. The latter is given by
\begin{equation} \label{parallel_F3}
\hat D_{\delta}F_{,\alpha\beta\gamma} = F_{,\alpha\beta\gamma\delta} - \frac12F^{\rho\sigma}(F_{,\alpha\beta\rho}F_{,\gamma\sigma\delta} + F_{,\alpha\gamma\rho}F_{,\beta\sigma\delta} + F_{,\beta\gamma\rho}F_{,\alpha\sigma\delta}).
\end{equation}
The last item now follows from the product structure of the pseudo-Riemannian manifold ${\cal D}$ established in Theorem \ref{splitting}.
\end{proof}

{\corollary \label{surface_chars} Let $f: M \to \mathbb R^n$ be a smooth non-degenerate centro-affine hypersurface immersion and let the function $F: {\cal D} = M \times \mathbb R_{++} \to \mathbb R$ be defined by \eqref{F_def}. Then the immersion $f$
\begin{itemize}
\item is a quadric if and only if $F_{,\alpha\beta\gamma} + 2F_{,\alpha\beta}F_{,\gamma} + 2F_{,\alpha\gamma}F_{,\beta} + 2F_{,\beta\gamma}F_{,\alpha} + 4F_{,\alpha}F_{,\beta}F_{,\gamma} = 0$,
\item is an affine hypersphere with center in the origin if and only if $F_{,\alpha\beta\gamma}F^{,\beta\gamma} + 2nF_{,\alpha} = 0$,
\item has parallel cubic form with respect to the Levi-Civita connection of the centro-affine metric if and only if
\end{itemize}
\begin{equation} \label{quasi_lin_PDE}
F_{,\alpha\beta\gamma\delta} = \frac12F^{,\rho\sigma}(F_{,\alpha\beta\rho}F_{,\gamma\delta\sigma} + F_{,\alpha\gamma\rho}F_{,\beta\delta\sigma} + F_{,\alpha\delta\rho}F_{,\beta\gamma\sigma}).
\end{equation} }

\begin{proof}
A covariant tensor $W$ on $M$ is zero if and only if its pullback $\Pi^*W$ on ${\cal D}$ is zero. Now $f$ is a quadric if and only if the cubic form $C$ vanishes, it is an affine hypersphere if and only if the Tchebycheff form $tr_hC$ vanishes, and it has parallel cubic form if and only if $\hat\nabla C$ vanishes. The corollary is then a direct consequence of Lemma \ref{correspondence}.
\end{proof}

The last two items of the corollary characterize affine hyperspheres and immersions with parallel cubic form conveniently by PDEs on the function $F$. The nonlinear PDE characterizing affine hyperspheres can actually be integrated, yielding the second order PDE $\det F'' = const\cdot e^{-2nF}$ of Monge-Amp\`ere type. This PDE is the non-convex analog of \cite[eq.(4.1), p.359]{ChengYau82}, which characterizes complete hyperbolic affine hyperspheres.

Let us introduce the difference tensor $K = D - \hat D$ the between the flat affine connection and the Levi-Civita connection of $F''$ on ${\cal D}$. By \eqref{Christoffel} we have
\begin{equation} \label{K_def}
K_{\alpha\beta}^{\gamma} = -\frac12F_{,\alpha\beta\delta}F^{,\gamma\delta}.
\end{equation}

{\lemma \label{K_parallel} Let $f: M \to \mathbb R^n$ be a smooth non-degenerate centro-affine hypersurface immersion. Then $f$ has parallel cubic form if and only if the difference tensor $K$ is parallel with respect to $\hat D$. }

\begin{proof}
By \eqref{parallel_F3} the derivative $F''' = D^3F$ is parallel with respect to $\hat D$ if and only if \eqref{quasi_lin_PDE} holds. Since the pseudo-metric $F''$ is parallel, we obtain that $K$ is parallel if and only if \eqref{quasi_lin_PDE} holds. Application of Corollary \ref{surface_chars} completes the proof.
\end{proof}

The situation is conceptually somewhat simpler if the centro-affine hypersurface is given implicitly as an integral manifold of a distribution rather than parametrically by a map $f: M \to \mathbb R^n$.

{\lemma \label{surf_chars_zeta} Let $Y \subset \mathbb R^n$ be an open subset, let $\zeta$ be a closed form on $Y$ such that $\zeta(x) \equiv 1$, where $x$ is the position vector field on $Y$, and assume that $D\zeta$ is non-degenerate. Denote by $\Psi$ the symmetric contravariant second order tensor which is the inverse of $D\zeta$. Let $\Delta$ be the involutive distribution on $Y$ given by the kernel of $\zeta$, and let $M \subset Y$ be an integral hypersurface of $\Delta$.

Then $\zeta$ satisfies the relation
\begin{equation} \label{FderXa}
(D\zeta)_{\alpha\beta}x^{\beta} = -\zeta_{\alpha},
\end{equation}
the Christoffel symbols of the Hessian pseudo-metric $D\zeta$ are given by $\Gamma_{\alpha\beta}^{\gamma} = \frac12\Psi^{\gamma\delta}(D^2\zeta)_{\alpha\beta\delta}$, and the difference tensor $K = D - \hat D$ between the flat affine connection $D$ on $\mathbb R^n$ and the Levi-Civita connection $\hat D$ of the pseudo-metric $D\zeta$ is given by
\begin{equation} \label{K_zeta}
K_{\alpha\beta}^{\gamma} = -\frac12\Psi^{\gamma\delta}(D^2\zeta)_{\alpha\beta\delta}.
\end{equation}
Moreover, $M$ is a non-degenerate centro-affine hypersurface with centro-affine metric given by the restriction of $D\zeta$ to $M$, and $M$
\begin{itemize}
\item is a quadric if $(D^2\zeta)_{\alpha\beta\gamma} + 2(D\zeta)_{\alpha\beta}\zeta_{\gamma} + 2(D\zeta)_{\alpha\gamma}\zeta_{\beta} + 2(D\zeta)_{\beta\gamma}\zeta_{\alpha} + 4\zeta_{\alpha}\zeta_{\beta}\zeta_{\gamma} = 0$,
\item is an affine hypersphere with center in the origin if $(D^2\zeta)_{\alpha\beta\gamma}\Psi^{\beta\gamma} + 2n\zeta_{\alpha} = 0$,
\item has parallel cubic form if
\end{itemize}
\[ (D^3\zeta)_{\alpha\beta\gamma\delta} = \frac12\Psi^{\rho\sigma}((D^2\zeta)_{\alpha\beta\rho}(D^2\zeta)_{\gamma\delta\sigma} + (D^2\zeta)_{\alpha\gamma\rho}(D^2\zeta)_{\beta\delta\sigma} + (D^2\zeta)_{\alpha\delta\rho}(D^2\zeta)_{\beta\gamma\sigma}).
\] }

\begin{proof}
First note that by $\zeta(x) = 1$ the distribution $\Delta$ is transversal to the position vector field and hence $M$ is centro-affine.

Equip $M$ with the topology of an immersed manifold. Let $y \in M$ be a point and $V \subset M$ a small enough neighbourhood of $y$ in $M$ with the following properties. Each ray in $\mathbb R^n$ intersects $V$ at most once, there exists $\varepsilon > 0$ such that $U = \bigcup_{|\lambda-1| < \varepsilon} \lambda V$ is contained in $Y$, and there exists a local potential $\Phi: U \to \mathbb R$ of $\zeta$ such that $\Phi|_V = 0$.

Consider $M$ as a centro-affine hypersurface immersion and let $\tilde f,F$ be the maps defined in \eqref{tilde_f},\eqref{F_def}, respectively. Define the set $U' = V \times (1-\varepsilon,1+\varepsilon) \subset {\cal D}$. Then $\tilde f[U'] = U$, and the restriction $\tilde f|_{U'}$ is injective.

We have $\zeta = D\Phi$ and $D_x\Phi = 1$. Integrating the latter relation along the rays of $\mathbb R^n$, we obtain the logarithmic homogeneity condition $\Phi(\alpha z) = \log\alpha + \Phi(z)$, whenever $z,\alpha z \in U$, $\alpha > 0$. Hence the pullback of $\Phi$ from $U$ to $U'$ by means of the restriction $\tilde f|_{U'}$ coincides with the function $F$, and the $k$-th derivative of $F$ on $U'$ is the pullback of the $(k-1)$-th derivative $D^{k-1}\zeta$ on $U$.

The claims of the lemma now follow from the second relation in Lemma \ref{FderX}, from equations \eqref{Christoffel} and \eqref{K_def}, and from Theorem \ref{splitting} and Corollary \ref{surface_chars}.
\end{proof}

The advantage of this description is that we can work directly on $\mathbb R^n$ with the 1-form $\zeta$ instead on ${\cal D}$ with the function $F$. For local considerations both approaches are equivalent. The class of centro-affine immersions admitting a global description as in Lemma \ref{surf_chars_zeta} is limited, however, for instance self-intersections are not allowed. We will see later that centro-affine immersions with parallel cubic form always allow a global description as an integral manifold of some involutive distribution.

\section{Jordan algebras} \label{sec_Jordan}

In this section we provide the necessary background on Jordan algebras. Most of the material in this section is taken from \cite{Koecher99}. Other references on Jordan algebras are \cite{Jacobson68} or \cite{McCrimmon}.

{\definition \label{Jordan_def} A {\sl Jordan algebra} $J$ over a field $\mathbb K$ is a vector space over $\mathbb K$ endowed with a bilinear operation $\bullet: J \times J \to J$ satisfying the following conditions:

i) commutativity: $x \bullet y = y \bullet x$ for all $x,y \in J$,

ii) Jordan identity: $x \bullet (x^2 \bullet y) = x^2 \bullet (x \bullet y)$ for all $x,y \in J$, where $x^2 = x \bullet x$.
}

Throughout the paper we assume that $\mathbb K$ is $\mathbb R$ or $\mathbb C$ and that $J$ is finite-dimensional.

\medskip

Let us denote the operator of multiplication with the element $x$ by $L_x$, $L_xy = x \bullet y = L_yx$. On $J$ we can define a linear form $t$ by
\begin{equation} \label{trace_form}
t(x) = tr\,L_x
\end{equation}
and a symmetric bilinear form $g$ by \cite[p.59]{Koecher99}
\begin{equation} \label{bilinear}
g(x,y) = t(x \bullet y) = tr\,L_{x \bullet y}.
\end{equation}
The form $g$ satisfies the relation \cite[Lemma III.4, p.59]{Koecher99}
\begin{equation} \label{3symmetry}
g(u \bullet v,w) = g(u,v \bullet w)
\end{equation}
for all $u,v,w \in J$. Equivalently, the operator $L_v$ is self-adjoint with respect to $g$ for all $v \in J$.

{\definition \cite[p.56]{Koecher99} The {\sl center} $Z(J)$ of a Jordan algebra $J$ is the set of all elements $z \in J$ such that $L_zL_x = L_xL_z$ for all $x \in J$. }

{\definition \cite[p.60]{Koecher99} A Jordan algebra $J$ is called {\sl semi-simple} if the bilinear form $g$ defined in \eqref{bilinear} is non-degenerate. }

{\definition \cite[p.64]{Koecher99} A Jordan algebra $J$ is called {\sl direct sum} of the subalgebras $J_1,\dots,J_r$, $J = \oplus_{k = 1}^r J_k$, if $J$ is the sum of $J_1,\dots,J_r$ as a vector space and $x \bullet y = 0$ for all $x \in J_k$, $y \in J_l$ with $k \not= l$.

A Jordan algebra $J$ is called {\sl simple} if it is semi-simple and cannot be represented as a nontrivial direct sum of subalgebras. }

Clearly the summands in a direct sum decomposition are ideals, i.e., $x \bullet y \in J_k$ for all $x \in J_k$, $y \in J$.

{\theorem \cite[Theorem 11, p.65]{Koecher99} \label{semi_simple_decomp} Let $J$ be a semi-simple Jordan algebra. Then there exist simple subalgebras $J_1,\dots,J_r \subset J$ such that $J = \oplus_{k=1}^r J_k$. Any two decompositions of $J$ into a direct sum of simple subalgebras are equal up to permutation of the summands. }

{\definition \cite[p.206]{Jacobson68} A simple Jordan algebra $J$ is called {\sl central simple} if the dimension of its center equals 1. }

{\lemma \cite[p.206]{Jacobson68} \label{lem_central_simple} Let $J$ be a simple Jordan algebra. Then its center $Z(J)$ is a field, and $J$ is a central-simple Jordan algebra if considered as an algebra over $Z(J)$. In particular, if $J$ is a real simple Jordan algebra, then it is either central-simple, or it is isomorphic to a central-simple Jordan algebra over $\mathbb C$. }

{\definition \cite{JvNW34} A Jordan algebra $J$ is called {\sl formally real} or {\sl Euclidean} if for all $x_1,\dots,x_r \in J$ the relation $\sum_{k=1}^r x_k^2 = 0$ implies $x_k = 0$, $k = 1,\dots,r$. }

{\theorem \cite[Theorem VI.12, p.118]{Koecher99} \label{th_posr} Let $J$ be a real Jordan algebra. Then the following are equivalent:

i) $J$ is formally real,

ii) the bilinear form $g$ defined by \eqref{bilinear} is positive definite,

iii) there exists a positive definite symmetric bilinear form $\sigma$ such that $\sigma(u \bullet v,w) = \sigma(u,v \bullet w)$ for all $u,v,w \in J$. }

Introduce the operator $P_x = 2L_x^2 - L_{x^2}$, which is quadratic in the parameter $x$. This operator satisfies the fundamental formula \cite[Theorem IV.1, p.73]{Koecher99}
\begin{equation} \label{fund_form}
P_{P_xy} = P_xP_yP_x
\end{equation}
for all $x,y \in J$.

An element $e \in J$ is called {\sl unit element} if $x \bullet e = x$ for all $x \in J$. A Jordan algebra possessing a unit element is called {\sl unital}. In a unital Jordan algebra $J$, $y = x^{-1}$ is called the {\sl inverse} of $x$ if $x \bullet y = e$ and $L_x,L_y$ commute \cite[p.67]{Koecher99}. If it exists, the inverse is unique \cite[Lemma III.5, p.66]{Koecher99} and satisfies $(x^{-1})^{-1} = x$ \cite[p.67]{Koecher99}. In this case we shall call $x$ {\sl invertible}. We have the following characterization of the invertible elements of $J$.

{\theorem \cite[Theorem III.12, p.67]{Koecher99} \label{th_invertible} Let $J$ be a unital Jordan algebra. An element $x \in J$ is invertible if and only if $\det P_x \not= 0$. In this case the inverse is given by $x^{-1} = P_x^{-1}x$, and
\[ P_{x^{-1}} = P_x^{-1},\qquad L_{x^{-1}} = L_xP^{-1}_x = P^{-1}_xL_x.
\] }

The derivative of the inverse is given by \cite[eq.(1), p.73]{Koecher99}
\begin{equation} \label{inv_der}
D_ux^{-1} = -P_x^{-1}u.
\end{equation}
Here $D_u$ denotes the derivative with respect to $x$ in the direction of $u$. We have \cite[Lemma IV.1, p.79]{Koecher99}
\begin{equation} \label{der_log_det}
D_u(\log\det P_x) = 2g(x^{-1},u)
\end{equation}
for all $u \in J$ and all invertible $x \in J$.

{\theorem \cite[Theorem III.9, p.63]{Koecher99} Every semi-simple Jordan algebra possesses a unit element. }

{\theorem \cite[Theorem III.10, p.64; pp.71--72]{Koecher99} \label{form_exist} Let $J$ be a unital Jordan algebra, and suppose that there exists a non-degenerate symmetric bilinear form $\gamma(u,v)$ on $J$ such that
\begin{equation} \label{gamma_sym}
\gamma(u \bullet v,w) = \gamma(u,v \bullet w)
\end{equation}
for all $u,v,w \in J$. Then the following are equivalent:

i) $\sigma(u,v)$ is a symmetric bilinear form on $J$ satisfying $\sigma(u \bullet v,w) = \sigma(u,v \bullet w)$ for all $u,v,w \in J$,

ii) there exists a central element $z \in Z(J)$ such that $\sigma(x,y) = \gamma(z \bullet x,y)$. }

In matrix form the relation in item ii) can be written as $\sigma = \gamma L_z$, or $L_z = \gamma^{-1}\sigma$. Note that in a unital algebra we have $L_x \not= 0$ for nonzero $x \in J$. Hence the map $x \mapsto L_x$ is injective and the element $z$ in item ii) is uniquely determined. For central elements we have $P_z = L_z^2$ and $z$ is invertible if and only if $L_z$ is. It follows that $\sigma$ is non-degenerate if and only if $z$ is invertible \cite[item (v), pp.71--72]{Koecher99}\footnote{The claim in \cite[item (v), pp.71--72]{Koecher99} that $J$ must be semi-simple is false. A counter-example will be given in Subsection \ref{subs_non_semi_simple}.}.

{\remark A symmetric bilinear form $\gamma$ satisfying \eqref{gamma_sym} is called a {\sl trace form} \cite[p.24]{Schaefer}. }

Let $J$ be a unital Jordan algebra. For $u \in J$, define a new multiplication $\bullet_u: J \times J \to J$ by \cite[p.76]{Koecher99}
\[ x \bullet_u y = x \bullet (y \bullet u) + y \bullet (x \bullet u) - (x \bullet y) \bullet u.
\]
Clearly this operation is bilinear and commutative, and the linear operator $L^{(u)}_x$ of multiplication with $x$ is given by
\[ L^{(u)}_x = L_xL_u - L_uL_x + L_{x \bullet u}.
\]
Define the corresponding quadratic operator $P^{(u)}_x = 2(L^{(u)}_x)^2 - L^{(u)}_{x \bullet_u x}$. 
Denote the algebra obtained by equipping $J$ with the multiplication $\bullet_u$ by $J^{(u)}$.

{\theorem \cite[pp.77--78]{Koecher99} Let $J$ be a unital Jordan algebra. For every $u \in J$ we have
\begin{equation} \label{P_prod}
P^{(u)}_x = P_xP_u, 
\end{equation}
and $J^{(u)}$ is a Jordan algebra. Moreover, if $J$ is semi-simple, then $J^{(u)}$ is semi-simple if and only if $u$ is invertible. }

If $u$ is invertible, then $[L_u,L_{u^{-1}}] = 0$ and $u^{-1}$ is the unit element of $J^{(u)}$.

{\definition \cite[p.57]{Jacobson68} For invertible $u \in J$ we call the Jordan algebra $J^{(u)}$ the {\sl $u$-isotope} of $J$ (in \cite{Koecher99} $J^{(u)}$ is called {\sl mutation}). }

\subsection{The set of invertible elements} \label{sec_inv}

Let $J$ be a unital real Jordan algebra. Note that $P_e$ is the identity matrix, and hence $\det P_e = 1$. Since $\det P_x$ is a polynomial in $x$, it follows by Theorem \ref{th_invertible} that the set ${\cal X}$ of invertible elements is open and dense in $J$. In this subsection we shall investigate this set and its symmetry properties.

Let $\Pi \subset GL(J)$ be the group generated by the transformations $P_u$, where $u$ varies in a small neighbourhood of the unit element. By \eqref{fund_form} the group $\Pi$ preserves the set ${\cal X}$ of invertible elements of $J$. Since $\Pi$ is connected, it preserves even every connected component of ${\cal X}$.

We set $u^0 = e$ and $u^{k+1} = L_uu^k$ recursively for all $u \in J$. Then we can define the {\sl exponential} $\exp(u) = \sum_{k=0}^{\infty} \frac{1}{k!} u^k = \exp(L_u)e$, which bijectively maps a neighbourhood of zero in $J$ to a neighbourhood of $e$ \cite[pp.82--83]{Koecher99}.

{\lemma \cite[Lemma IV.4, p.83]{Koecher99} \label{LexpP} For every $u \in J$ we have $P_{\exp(u)} = \exp(2L_u)$. }

From this we have the following result.

{\lemma \label{lem_gen} The group $\Pi$ is generated by the 1-dimensional subgroups $\exp(tL_w)$, $w \in J$. \qed }

The action of the subgroup $\exp(tL_w)$ on $J$ generates a flow with tangent vector field $X_w(x) = L_wx = x \bullet w = L_xw$. If $L_x$ is non-degenerate, then the vectors $\{ X_w(x) \,|\, w \in J \}$ span the whole tangent space $T_xJ$, and the orbit of $x$ under the action of $\Pi$ has full dimension.

{\lemma \cite[Lemma V.6, p.106]{Koecher99} Let $J$ be a real Jordan algebra with unit element $e$ and let $x \in J$ be invertible. Then in any neighbourhood of $e$ there exists an element $w$ such that $\det L_{P_wx} \not= 0$. }

It follows that the orbit of $x$ has full dimension for all invertible $x$. By a topological argument \cite[pp.110--111]{Koecher99} we get the following result.

{\theorem \cite[Theorem VI.2, p.110]{Koecher99} \label{transitive_action} Let $J$ be a unital real Jordan algebra. Then the group $\Pi$ acts transitively on every connected component of the set ${\cal X}$ of invertible elements of $J$. }

In the context of \cite[Chapter VI]{Koecher99} this theorem was proven for semi-simple $J$, but as we have sketched above, the proof is valid for every unital $J$.

\medskip

We now consider the case when the Jordan algebra $J$ is semi-simple. First we introduce the notion of an {\sl $\omega$-domain}.

Let $Y$ be a nonempty open connected subset of an $n$-dimensional real vector space $V$ such that $\lambda y \in Y$ for all $\lambda > 0$ and $y \in Y$. Let $\omega$ be a continuous real-valued function defined on the closure of $Y$ which is analytic and positive on $Y$, vanishes on the boundary $\partial Y$, satisfies $\omega(\lambda y) = \lambda^n\omega(y)$ for all $\lambda > 0$, $y \in Y$, and such that the Hessian of $\log\omega$ is non-degenerate on $Y$.

Then $(-\log\omega)''$ defines a Hessian pseudo-metric $\sigma$ on $Y$. Since $Y$ is a subset of the vector space $V$, $\sigma$ can be seen as a map from $Y$ to the space of symmetric bilinear forms on $V$. Let us note $\sigma_y$ for the value of this map at $y \in Y$. For points $x,y \in Y$, the map $H^x_y$ given by the matrix $\sigma_x^{-1}\sigma_y$ is an element of $GL(V)$. We shall assume that actually $H^x_y \in \Sigma(Y,\omega)$ for every $x,y \in Y$, where $\Sigma(Y,\omega)$ is the group of linear transformations $W \in GL(V)$ such that $WY = Y$ and $\omega(Wx) = |\det W|\cdot\omega(x)$ for all $x \in Y$.

{\definition \cite[pp.35--36]{Koecher99} A pair $(Y,\omega)$ satisfying the above assumptions is called {\sl $\omega$-domain}. }

Note that for $W \in \Sigma(Y,\omega)$ we have $-\log\omega(Wx) = -\log|\det W|-\log\omega(x)$, and therefore $\Sigma(Y,\omega)$ preserves the pseudo-metric $\sigma = (-\log\omega)''$.

{\theorem \label{domain_algebra} Assume the above notations. Let $J$ be a real semi-simple Jordan algebra and ${\cal X}$ the set of its invertible elements. Then the connected components of ${\cal X}$ are {\sl $\omega$-domains} \cite[Corollary VI.4, p.112]{Koecher99}. The function $\omega$ of the connected component ${\cal Y}$ of $e$ in ${\cal X}$ is given by $\omega(y) = \sqrt{|\det P_y|}$ \cite[eq. VI.3, p.114]{Koecher99}, and the maps $H^e_y$ by $P^{-1}_y$. The pseudo-metric $\sigma$ at the unit element is given by the bilinear form \eqref{bilinear}, $\sigma_e = g$ \cite[p.114]{Koecher99}.


On the other hand, for every $\omega$-domain $Y \subset V$, where $V$ is a real vector space, and every point $c \in Y$ there exists a real semi-simple Jordan algebra $J$ on $V$ with unit element $e = c$ such that $Y$ is the connected component ${\cal Y}$ of $e$ in ${\cal X}$ \cite[Theorem VI.6, p.114]{Koecher99}. }

In particular, the set ${\cal Y}$ inherits the Hessian pseudo-metric of the $\omega$-domain. By \cite[p.79; Theorem VI.1, p.110]{Koecher99} the group $\Pi$ is a subgroup of $\Sigma({\cal Y},\omega)$ and hence consists of isometries. Thus ${\cal Y}$ is a homogeneous space. It is actually even a symmetric space, because the map $x \mapsto x^{-1}$ is an involution of ${\cal Y}$ \cite[Theorem 1]{Koecher70}.

Finally we consider the other connected components of ${\cal X}$. Let ${\cal Y}' \subset {\cal X}$ be such a component, and let $v \in {\cal Y}'$. Then $v$ is invertible. Let $u$ be its inverse and consider the isotope $J^{(u)}$. By \eqref{P_prod} the sets of invertible elements in $J$ and $J^{(u)}$ coincide. Since $v = u^{-1}$ is the unit element in $J^{(u)}$, the component ${\cal Y}'$ is then the connected component in ${\cal X}$ of the unit element in $J^{(u)}$. By \eqref{P_prod} its $\omega$-function is just a multiple of the original $\omega$-function.

%
%
%
%
%

\section{Immersions and Jordan algebras} \label{sec_correspondence}

We are now in a position to establish a connection between non-degenerate centro-affine hypersurface immersions with parallel cubic form and real unital Jordan algebras with non-degenerate trace form.

\subsection{Jordan algebras defined by immersions}

In this subsection we show that a non-degenerate centro-affine hypersurface immersion $f: M \to \mathbb R^n$ with parallel cubic form together with a point $y \in {\cal D} = M \times \mathbb R_{++}$ equip the tangent space $T_y{\cal D}$ with the structure of a unital Jordan algebra $J$.

Let $f: M \to \mathbb R^n$ be a connected non-degenerate centro-affine hypersurface immersion and let the function $F: {\cal D} \to \mathbb R$ be defined by \eqref{F_def}. By Corollary \ref{surface_chars} $f$ has parallel cubic form if and only if $F$ obeys the quasi-linear fourth-order PDE \eqref{quasi_lin_PDE}.

Let us deduce the integrability condition of this PDE. Introduce local affine coordinates $x^{\alpha}$ on ${\cal D}$. Differentiating \eqref{quasi_lin_PDE} with respect to the coordinate $x^{\eta}$ and substituting the appearing fourth order derivatives of $F$ by the right-hand side of \eqref{quasi_lin_PDE}, we obtain after simplification
\begin{eqnarray*}
F_{,\alpha\beta\gamma\delta\eta} &=& \frac14F^{,\rho\sigma}F^{,\mu\nu}\left(F_{,\beta\eta\nu}F_{,\alpha\rho\mu}F_{,\gamma\delta\sigma} + F_{,\alpha\eta\mu}F_{,\rho\beta\nu}F_{,\gamma\delta\sigma} + F_{,\gamma\eta\nu}F_{,\alpha\rho\mu}F_{,\beta\delta\sigma} + F_{,\alpha\eta\mu}F_{,\rho\gamma\nu}F_{,\beta\delta\sigma} \right. \\
&& + F_{,\beta\eta\nu}F_{,\gamma\rho\mu}F_{,\alpha\delta\sigma} + F_{,\gamma\eta\mu}F_{,\rho\beta\nu}F_{,\alpha\delta\sigma} + F_{,\beta\eta\nu}F_{,\delta\rho\mu}F_{,\alpha\gamma\sigma} + F_{,\delta\eta\mu}F_{,\rho\beta\nu}F_{,\alpha\gamma\sigma} \\
&& \left. + F_{,\delta\eta\nu}F_{,\alpha\rho\mu}F_{,\beta\gamma\sigma} + F_{,\alpha\eta\mu}F_{,\rho\delta\nu}F_{,\beta\gamma\sigma} + F_{,\delta\eta\nu}F_{,\gamma\rho\mu}F_{,\alpha\beta\sigma} + F_{,\gamma\eta\mu}F_{,\rho\delta\nu}F_{,\alpha\beta\sigma}\right).
\end{eqnarray*}
The right-hand side must be symmetric in all 5 indices. Commuting the indices $\delta,\eta$ and equating the resulting expression with the original one we obtain
\begin{eqnarray*}
\lefteqn{F^{,\rho\sigma}F^{,\mu\nu}\left(F_{,\beta\eta\nu}F_{,\delta\rho\mu}F_{,\alpha\gamma\sigma} + F_{,\alpha\eta\mu}F_{,\rho\delta\nu}F_{,\beta\gamma\sigma} + F_{,\gamma\eta\mu}F_{,\rho\delta\nu}F_{,\alpha\beta\sigma} \right.} \\
&& \left. - F_{,\beta\delta\nu}F_{,\eta\rho\mu}F_{,\alpha\gamma\sigma} - F_{,\alpha\delta\mu}F_{,\rho\eta\nu}F_{,\beta\gamma\sigma} - F_{,\gamma\delta\mu}F_{,\rho\eta\nu}F_{,\alpha\beta\sigma}\right) = 0.
\end{eqnarray*}
Raising the index $\eta$, we get by virtue of \eqref{K_def} the integrability condition
\[ K_{\alpha\mu}^{\eta}K_{\delta\rho}^{\mu}K_{\beta\gamma}^{\rho} + K_{\beta\mu}^{\eta}K_{\delta\rho}^{\mu}K_{\alpha\gamma}^{\rho} + K_{\gamma\mu}^{\eta}K_{\delta\rho}^{\mu}K_{\alpha\beta}^{\rho} = K_{\alpha\delta}^{\mu}K_{\rho\mu}^{\eta}K_{\beta\gamma}^{\rho} + K_{\beta\delta}^{\mu}K_{\rho\mu}^{\eta}K_{\alpha\gamma}^{\rho} + K_{\gamma\delta}^{\mu}K_{\rho\mu}^{\eta}K_{\alpha\beta}^{\rho}.
\]
This condition is satisfied if and only if
\[ K_{\alpha\mu}^{\eta}K_{\delta\rho}^{\mu}K_{\beta\gamma}^{\rho}u^{\alpha}u^{\beta}u^{\gamma}v^{\delta} = K_{\alpha\delta}^{\mu}K_{\rho\mu}^{\eta}K_{\beta\gamma}^{\rho}u^{\alpha}u^{\beta}u^{\gamma}v^{\delta}
\]
for all tangent vector fields $u,v$ on ${\cal D}$, which can be written as
\begin{equation} \label{int_cond}
K(K(K(u,u),v),u) = K(K(u,v),K(u,u)).
\end{equation}

{\theorem \label{main1} Let $f: M \to \mathbb R^n$ be a non-degenerate centro-affine hypersurface immersion with parallel cubic form and let the function $F: {\cal D} = M \times \mathbb R_{++} \to \mathbb R$ be defined by \eqref{F_def}. Let $y \in {\cal D}$ be a point and let $\bullet: T_y{\cal D} \times T_y{\cal D} \to T_y{\cal D}$ be the multiplication $(u,v) \mapsto K(u,v)$ defined by the tensor $K$ at $y$. Let $\gamma$ be the symmetric non-degenerate bilinear form defined on $T_y{\cal D}$ by the pseudo-metric $F''$.

Then the tangent space $T_y{\cal D}$, equipped with the multiplication $\bullet$, is a unital real Jordan algebra $J$, and the position vector $e$ is its unit element. The form $\gamma$ is a trace form and $\gamma(e,e) = -1$. There exists a unique central element $z \in Z(J)$ such that $g(u,v) = \gamma(z \bullet u,v)$ for all $u,v \in J$, where $g$ is the bilinear form defined by \eqref{bilinear}. }

\begin{proof}
Assume the conditions of the theorem. The tensor $K_{\alpha\beta}^{\gamma}$ is symmetric in the indices $\alpha,\beta$, hence the multiplication $\bullet$ is commutative. With $u^2 = u \bullet u$ condition \eqref{int_cond} becomes equivalent to
\[ u \bullet (u^2 \bullet v) = u^2 \bullet (u \bullet v),
\]
which is the Jordan identity in Definition \ref{Jordan_def}. Thus $T_y{\cal D}$, equipped with the multiplication $\bullet$, is a Jordan algebra $J$.

By the third relation in Lemma \ref{FderX} we have
\begin{equation} \label{pos_identity}
(u \bullet e)^{\gamma} = K^{\gamma}_{\alpha\beta}u^{\alpha}e^{\beta} = -\frac12F_{,\alpha\beta\delta}F^{,\gamma\delta}u^{\alpha}e^{\beta} = F_{,\alpha\delta}F^{,\gamma\delta}u^{\alpha} = u^{\gamma},
\end{equation}
for all $u \in J$, and the position vector $e$ is the unit element of $J$.

For any vectors $u,v,w \in J$ we have
\begin{eqnarray} \label{gamma_symmetric}
\gamma(u \bullet v,w) &=& F_{,\beta\gamma}K_{\delta\rho}^{\beta}u^{\delta}v^{\rho}w^{\gamma} = -\frac12F_{,\beta\gamma}F_{,\delta\rho\sigma}F^{,\sigma\beta}u^{\delta}v^{\rho}w^{\gamma} = -\frac12F_{,\delta\rho\gamma}u^{\delta}v^{\rho}w^{\gamma} \nonumber \\ &=& -\frac12F_{,\beta\delta}u^{\delta}F_{,\rho\gamma\sigma}F^{,\sigma\beta}v^{\rho}w^{\gamma} =
F_{,\delta\beta}u^{\delta}K_{\rho\gamma}^{\beta}v^{\rho}w^{\gamma} = \gamma(u,v \bullet w).
\end{eqnarray}
Hence the form $\gamma$ satisfies \eqref{gamma_sym} and verifies the conditions of Theorem \ref{form_exist}.

By \eqref{3symmetry} the bilinear form $g$ verifies the assumptions on the form $\sigma$ in Theorem \ref{form_exist}, which implies existence and uniqueness of a central element $z$ such that $g(u,v) = \gamma(z \bullet u,v)$ for all $u,v$.

Finally, the relation $\gamma(e,e) = -1$ follows from the fifth relation in Lemma \ref{FderX}.
\end{proof}

Let us compute the central element $z$ in terms of the derivatives of $F$. We have $tr\,L_v = g(e,v) = \gamma(z \bullet e,v) = \gamma(z,v)$ for all $v \in J$. This yields $K_{\alpha\gamma}^{\gamma} = F_{,\alpha\gamma}z^{\gamma}$, and hence
\begin{equation} \label{z_expr}
z^{\delta} = -\frac12F_{,\alpha\beta\gamma}F^{,\beta\gamma}F^{,\alpha\delta}.
\end{equation}

{\lemma \label{AHS} Assume the conditions of Theorem \ref{main1}. Then $f$ is an affine hypersphere with center in the origin if and only if $z = -ne$, or equivalently, $g = -n\gamma$. In this case $J$ is semi-simple. }

\begin{proof}
From \eqref{z_expr} and the first relation in \eqref{FderX2} we have
\[ z^{\delta} = -\frac12(\Pi^*T_{\alpha} - 2nF_{,\alpha})F^{,\alpha\delta} = -\frac12F^{,\alpha\delta}\Pi^*T_{\alpha} - ne^{\delta},
\]
where $\Pi^*T_{\alpha}$ is the pullback of the Tchebycheff form from Lemma \ref{correspondence}.

If $f$ is an affine hypersphere, then by Corollary \ref{surface_chars} we have $\Pi^*T_{\alpha} = 0$ and hence $z = -ne$.

If, on the other hand, $z = -ne$, then $\Pi^*T = 0$ at $y = (\xi,\lambda) \in {\cal D}$. Hence the Tchebycheff form $T$ vanishes at $\xi \in M$, and the cubic form $C$ of $M$ is traceless at the point $\xi$. But $C$ is parallel, and hence must be traceless everywhere on $M$. Thus $M$ is an affine hypersphere with center in the origin.

Finally, the relation $z = -ne$ implies that $L_z = -nI$ is invertible. Hence $g = \gamma L_z$ is non-degenerate and $J$ is semi-simple.
\end{proof}

{\corollary Assume the conditions of Theorem \ref{main1}. If the Jordan algebra $J$ defined by the hypersurface immersion $f$ is central-simple, then $M$ is an affine hypersphere with center in the origin. }

\begin{proof}
If $J$ is central-simple, then $z = \mu e$ for some constant $\mu$. Since $g(e,e) = t(e) = n$, we have $n = \gamma(z \bullet e,e) = \mu \gamma(e,e) = -\mu$. The claim now follows from Lemma \ref{AHS}.
\end{proof}

{\lemma Assume the conditions of Theorem \ref{main1}. If the form $\gamma$ is negative definite, i.e., $f$ is locally strongly convex of hyperbolic type, then the Jordan algebra $J$ is formally real. }

\begin{proof}
If $\gamma$ is negative definite, then $\sigma = -\gamma$ is positive definite and $J$ is formally real by Theorem \ref{th_posr}.
\end{proof}

We shall now provide the equivalent of Theorem \ref{main1} in the case when the centro-affine hypersurface is given as an integral hypersurface of an involutive distribution.

{\theorem \label{main1b} Assume the conditions and notations of Lemma \ref{surf_chars_zeta}. Suppose that $M$ has parallel cubic form as a centro-affine hypersurface. Let $e \in M$ be a point and let $\bullet: \mathbb R^n \times \mathbb R^n \to \mathbb R^n$ be the multiplication $(u,v) \mapsto K(u,v)$ defined by the difference tensor $K$ at $e$. Let further $I_e: T_eY \to \mathbb R^n$ be the canonical isomorphism between $\mathbb R^n$ and the tangent space at $e$. Define a symmetric bilinear form $\gamma$ on $\mathbb R^n$ by the pullback $(I_e^{-1})_*(D\zeta)$.

Then $\mathbb R^n$, equipped with the multiplication $\bullet$, is a unital real Jordan algebra $J$, and the vector $e$ is its unit element. The form $\gamma$ is a trace form and $\gamma(e,e) = -1$. There exists a unique central element $z \in Z(J)$ such that $g(u,v) = \gamma(z \bullet u,v)$ for all $u,v \in J$, where $g$ is the bilinear form defined by \eqref{bilinear}. }

\begin{proof}
As in the proof of Lemma \ref{surf_chars_zeta}, we find a neighbourhood $U$ of $e$ in $\mathbb R^n$ and a neighbourhood $U'$ of $a = (e,1)$ in ${\cal D}$ such that $\tilde f[U'] = U$, the restriction $\tilde f|_{U'}$ is injective, and there exists a local potential $\Phi$ of $\zeta$ on $U$ such that the restriction $F|_{U'}$ is the pullback of $\Phi$ from $U$ by means of $\tilde f|_{U'}$.

Then $\tilde f(a) = e$, $\tilde f_*[T_a{\cal D}] = T_eY$, and $I_e \circ \tilde f_*: T_a{\cal D} \to \mathbb R^n$ is an isomorphism that takes the position vector in $T_a{\cal D}$ to $e$. The claims of Theorem \ref{main1b} now follow by carrying over to $\mathbb R^n$ the structures $\bullet,\gamma$ defined on $T_a{\cal D}$ by virtue of Theorem \ref{main1}.
\end{proof}

\subsection{Immersions defined by Jordan algebras} \label{sec_imm_alg}

In this subsection we consider the opposite direction and show that every pair $(J,\gamma)$, where $J$ is a real unital Jordan algebra of dimension $n \geq 2$ and $\gamma$ is a non-degenerate trace form on $J$ satisfying $\gamma(e,e) = -1$, defines an invariant involutive distribution $\Delta$ on the connected component ${\cal Y}$ of $e$ in the set ${\cal X}$ of invertible elements, such that the integral hypersurfaces of $\Delta$ are centro-affine and have parallel cubic form.

For a pair $(J,\gamma)$ as above, let us define a 1-form $\zeta$ on ${\cal Y}$. At the point $x \in {\cal Y}$ we set
\begin{equation} \label{zeta_def}
\zeta(u) = -\gamma(u,x^{-1}).
\end{equation}
Let $\Gamma = \{ u \in J \,|\, \gamma(u,e) = 0 \}$ be the orthogonal complement of $e$ with respect to $\gamma$.

{\lemma \label{lem_zeta} The form $\zeta$ evaluates to 1 on the position vector field $x$, is closed, and invariant with respect to the action of the group $\Pi$. The kernel $\Delta$ of $\zeta$ is an involutive $(n-1)$-dimensional distribution on ${\cal Y}$. The derivative $D\zeta$ is non-degenerate. The form $\zeta$ annihilates the vector field $X_w(x) = x \bullet w$ for every $w \in \Gamma$. Both the distribution $\Delta$ and the form $D\zeta$ are invariant with respect to the map $x \mapsto x^{-1}$ on ${\cal Y}$. All assertions except the relation $\zeta(x) \equiv 1$ remain valid if we drop the condition $\gamma(e,e) = -1$. }

\begin{proof}
Suppose that $\zeta$ is zero at some point $x \in {\cal Y}$. Then $x^{-1} \not= 0$ is in the kernel of $\gamma$, which contradicts the non-degeneracy of $\gamma$. Hence $\zeta$ is nowhere zero and its kernel $\Delta$ is an $(n-1)$-dimensional distribution.

At $x \in {\cal Y}$ the derivative of $\zeta$ in the direction of a vector $v \in J$ is by virtue of \eqref{inv_der} given by
\begin{equation} \label{der_zeta}
(D_v\zeta)(u) = -\gamma(u,D_vx^{-1}) = \gamma(u,P_x^{-1}v).
\end{equation}
But the operator $P_x^{-1}$ is self-adjoint with respect to $\gamma$, because the operators $L_x,L_{x^2}$ are self-adjoint by virtue of \eqref{gamma_sym}. Thus $(D_v\zeta)(u) = (D_u\zeta)(v)$, $D\zeta$ is symmetric in its two arguments, and the exterior derivative of $\zeta$ vanishes. It follows that $\zeta$ is closed and the distribution $\Delta$ involutive.

By \eqref{der_zeta} the matrix of $D\zeta$ is given by $\gamma P_x^{-1}$. But both $\gamma$ and $P_x$ are non-degenerate. Hence $D\zeta$ is non-degenerate.

We now show that $\zeta$ is invariant with respect to the action of $\Pi$. By Lemma \ref{lem_gen} it is sufficient to show that for every $w \in J$ the Lie derivative of $\zeta$ in the direction of the vector field $X_w(x) = x \bullet w$ vanishes. We have
\begin{eqnarray*}
({\cal L}_{x \bullet w}\zeta)(u) &=& (D_{x \bullet w}\zeta)(u) + \zeta(D_u(x \bullet w)) = \gamma(u,P_x^{-1}(x \bullet w)) + \zeta(u \bullet w) \\ &=& \gamma(u,P_x^{-1}L_xw) - \gamma(u \bullet w,x^{-1}) = \gamma(u,L_{x^{-1}}w) - \gamma(u,w \bullet x^{-1}) = 0.
\end{eqnarray*}
Here we used \eqref{der_zeta} for the second equality and Theorem \ref{th_invertible} for the fourth.

Let $w \in \Gamma$. Then $\zeta(X_w) = -\gamma(x \bullet w,x^{-1}) = -\gamma(w,x \bullet x^{-1}) = -\gamma(w,e) = 0$ by the definition of $\Gamma$.

Consider vectors $u,v \in T_xJ$. By \eqref{inv_der} the images $\tilde u,\tilde v \in T_{x^{-1}}J$ of $u,v$ under the differential of the map $x \mapsto x^{-1}$ are given by $\tilde u = -P_x^{-1}u$, $\tilde v = -P_x^{-1}v$. We then have
\begin{eqnarray*}
\zeta(\tilde u) &=& -\gamma(\tilde u,x) = \gamma(P_x^{-1}u,x) = \gamma(u,P_x^{-1}x) = \gamma(u,x^{-1}) = -\zeta(u), \\
(D\zeta)(\tilde u,\tilde v) &=& \gamma(\tilde u,P_{x^{-1}}^{-1}\tilde v) = \gamma(P_x^{-1}u,P_xP_x^{-1}v) = \gamma(u,P_x^{-1}v) = (D\zeta)(u,v).
\end{eqnarray*}
Here in the third equalities of both lines we used that $P_x^{-1}$ is self-adjoint with respect to $\gamma$. In the first and last equalities of the second line we used \eqref{der_zeta}. In the second equality of the second line and in the fourth equality of the first line we used Theorem \ref{th_invertible}. From the first line it follows that $\Delta$ is invariant with respect to the non-linear transformation $x \mapsto x^{-1}$, while the second line implies that $D\zeta$ is invariant.

Up to now we did not use the condition $\gamma(e,e) = -1$. With this condition we have for $x \in {\cal Y}$ that $\zeta(x) = -\gamma(x,x^{-1}) = -\gamma(x \bullet e,x^{-1}) = -\gamma(e,x \bullet x^{-1}) = -\gamma(e,e) = 1$. This completes the proof.
\end{proof}

Define the subspace $L(\Gamma) = \{ L_w \,|\, w \in \Gamma \} \subset L(J) = \{ L_w \,|\, w \in J \}$ and let $\Theta \subset \Pi$ be the subgroup generated by the 1-dimensional subgroups $\exp(tL_w)$, $w \in \Gamma$. Let ${\cal L}_{\Theta},{\cal L}_{\Pi}$ be the Lie algebras generated by $L(\Gamma),L(J)$, respectively.

{\lemma \label{Gamma_z} The group $\Theta$ preserves the maximal integral manifolds of $\Delta$. The Lie algebra ${\cal L}_{\Pi}$ can be expressed as a direct sum ${\cal L}_{\Theta} \oplus (\mathbb R\cdot I)$, with $I$ being the identity matrix. }

\begin{proof}
By Lemma \ref{lem_zeta} the vector fields $X_w$ are tangent to $\Delta$ for $w \in \Gamma$. This implies the first claim of the lemma.

Now note that $e \not\in \Gamma$, because $\gamma(e,e) = -1 \not= 0$. Hence $J = \Gamma \oplus (\mathbb R\cdot e)$, which implies $L(J) = L(\Gamma) \oplus (\mathbb R\cdot I)$. It follows that $[L(J),L(J)] = [L(\Gamma),L(\Gamma)]$, $[L(J),[L(J),L(J)]] = [L(\Gamma),[L(\Gamma),L(\Gamma)]]$ etc. Hence ${\cal L}_{\Pi} = {\cal L}_{\Theta} + (\mathbb R\cdot I)$.

It rests to show that $I$ is not an element of ${\cal L}_{\Theta}$. The position vector field $Ix = x$ generated by the action of the subgroup $\exp(tI)$ is not tangent to the involutive distribution $\Delta$. However, all the generators of ${\cal L}_{\Theta}$ induce vector fields which are tangent to $\Delta$, and thus the vector fields $Ax$ must be tangent to $\Delta$ for all $A \in {\cal L}_{\Theta}$. This implies that $I$ does not belong to ${\cal L}_{\Theta}$.
\end{proof}

Let us define the subgroup $S\Pi = \Pi \cap SL(J)$ of elements with determinant 1. Clearly $S\Pi$ is connected. Since $P_{\alpha e} = \alpha^2 P_e = \alpha^2 I$ is an element of $\Pi$ for all $\alpha > 0$, and $\det G > 0$ for all $G \in \Pi$ by the connectedness of $\Pi$, we have that $\Pi$ is the direct product of $S\Pi$ with the central subgroup $\mathbb R_{++} I = \{ \alpha I \,|\, \alpha > 0 \}$. The group $S\Pi$ is generated by the 1-dimensional subgroups $\exp(tL_w)$ satisfying $tr\,L_w = t(w) = 0$. Let ${\cal L}_{S\Pi}$ be the Lie algebra generated by $\{ L_w \,|\, t(w) = 0 \}$. As in the proof of Lemma \ref{Gamma_z}, we have the decomposition ${\cal L}_{\Pi} = {\cal L}_{S\Pi} \oplus (\mathbb R\cdot I)$.

{\lemma \label{covering} Let $G \in S\Pi$ be an arbitrary element. Then there exist $\alpha > 0$ and $G' \in \Theta$ such that $G' = \alpha G$. }

\begin{proof}
Define the group homomorphism $\pi: \Theta \to S\Pi$ by $\pi(G) = (\det G)^{-1/n}G$. Its differential $d\pi|_{Id}: {\cal L}_{\Theta} \to {\cal L}_{S\Pi}$ at the identity element is given by $d\pi|_{Id}: A \mapsto A - (n^{-1}tr\,A)I$. We have ${\cal L}_{\Theta} \oplus (\mathbb R\cdot I) = {\cal L}_{S\Pi} \oplus (\mathbb R\cdot I)$, and $d\pi|_{Id}$ is an isomorphism of the Lie algebras ${\cal L}_{\Theta},{\cal L}_{S\Pi}$. By \cite[Proposition 3.26, p.100]{Warner} $\pi$ must then be a covering map. In particular, it is surjective. The claim of the lemma now easily follows.
\end{proof}

Define the function $\varphi(u) = \det P_u$ on $J$. Since $\varphi(P_uv) = \varphi^2(u)\varphi(v)$ \cite[p.78]{Koecher99}, the function $\varphi$ is invariant under the action of $S\Pi$. In particular, the level surfaces $\varphi_c = \{ x \in {\cal Y} \,|\, \varphi(x) = c \}$, $c > 0$, are the orbits of the action of $S\Pi$ on ${\cal Y}$. Let $\pi: y \mapsto (\varphi(y))^{-\frac{1}{2n}}y$ be the projection of ${\cal Y}$ onto the level surface $\varphi_1$ along the rays of ${\cal Y}$.

{\theorem \label{alg_to_imm} Let $M$ be a maximal integral manifold of the distribution $\Delta$. Then $M$ is a homogeneous, symmetric, non-degenerate, affine complete centro-affine hypersurface with parallel cubic form. Its conic hull $\bigcup_{\lambda > 0} \lambda M$ equals ${\cal Y}$. The restriction $\pi|_M$ is a covering map, and different sheets of the covering are related by homothety. }

\begin{proof}
That $M$ is centro-affine follows from transversality of $\Delta$ to the position vector field.

By Theorem \ref{transitive_action} the group $\Pi$ acts transitively on ${\cal Y}$. Let $y,y' \in M$ be arbitrary points. Then there exists an element $G \in \Pi$ such that $Gy = y'$. But by Lemma \ref{lem_zeta} $G$ preserves the form $\zeta$ and hence also the distribution $\Delta$. Thus $G$ takes the maximal integral manifold of $\Delta$ through $y$ to the maximal integral manifold through $y'$. But both manifolds coincide with $M$, and $G$ restricts to a diffeomorphism of $M$. Thus there exists a subgroup of $\Pi$ that acts transitively on $M$.

By Lemma \ref{surf_chars_zeta} the centro-affine metric on $M$ is given by the restriction of $D\zeta$ to $M$. By Lemma \ref{lem_zeta} $D\zeta$ is non-degenerate, hence so is the hypersurface $M$. Moreover, by Lemma \ref{lem_zeta} the form $\zeta$ and hence its derivatives are invariant with respect to the action of the group $\Pi$. Thus $\Pi$ consists of isometries of ${\cal Y}$ and $M$ is a homogeneous space. It also follows that $M$ is affine complete.

By Lemma \ref{lem_zeta} the mapping $x \mapsto x^{-1}$ preserves both the distribution $\Delta$ on ${\cal Y}$ and the pseudo-metric defined by $D\zeta$. Since $e$ is a fixed point of this map, the map is an involution of the maximal integral manifold $M'$ of $\Delta$ through $e$. Thus $M'$ is a symmetric space. But $M$ is a homothetic image of $M'$ and hence also a symmetric space.

Let now $y \in {\cal Y}$ be an arbitrary point. Choose a point $x \in M$. Then there exists $G \in \Pi$ such that $Gx = y$. We have $(\det G)^{-1/n}G \in S\Pi$. By Lemma \ref{covering} there exists $\alpha > 0$ and $G' \in \Theta$ such that $G' = \alpha(\det G)^{-1/n}G$. Since the action of $\Theta$ on ${\cal Y}$ preserves $M$, we have $G'x = \alpha(\det G)^{-1/n}Gx = \alpha(\det G)^{-1/n}y \in M$. Hence $y$ is in the conic hull of $M$ and this conic hull must be equal to ${\cal Y}$.

It follows that the projection $\pi|_M: M \to \varphi_1$ is surjective. Since the distribution $\Delta$ is invariant with respect to homotheties, every simply connected neighbourhood in $\varphi_1$ is evenly covered by $\pi|_M$. Thus $\pi|_M$ is a covering map, and the sheets over evenly covered neighbourhoods are related by homothety.

Let us show that $M$ has parallel cubic form. Let $u,v \in J$ be arbitrary vectors. Note that $D_vP_x = 2L_xL_v + 2L_vL_x - 2L_{x \bullet v}$. By virtue of \eqref{der_zeta} we then have at the point $x \in {\cal Y}$
\begin{eqnarray*}
(D\zeta)(u,v) &=& \gamma(u,P_x^{-1}v), \\
(D^2\zeta)(u,v,v) &=& -2\gamma(u,P_x^{-1}(L_xL_v + L_vL_x - L_{x \bullet v})P_x^{-1}v), \\
(D^3\zeta)(u,v,v,v) &=& 8\gamma(u,P_x^{-1}(L_xL_v + L_vL_x - L_{x \bullet v})P_x^{-1}(L_xL_v + L_vL_x - L_{x \bullet v})P_x^{-1}v) \\
&& -2\gamma(u,P_x^{-1}P_vP_x^{-1}v).
\end{eqnarray*}
Specifying to $x = e$, we obtain
\begin{eqnarray} \label{zeta_deriv_e}
(D\zeta)(u,v) &=& \gamma(u,v), \\
(D^2\zeta)(u,v,v) &=& -2\gamma(u,L_vv) = -2\gamma(u,v^2), \nonumber\\
(D^3\zeta)(u,v,v,v) &=& 8\gamma(u,L_v^2v) - 2\gamma(u,P_vv) = 6\gamma(u,v^3) = 6\gamma(v \bullet u,v^2). \nonumber
\end{eqnarray}
The matrix of $D\zeta$ at $e$ is then given by $\gamma$, and its inverse by $\gamma^{-1}$. The coordinate vector of the 1-form $(D^2\zeta)(\cdot,v,v)$ is given by $-2\gamma v^2$. If we consider the inverse $(D\zeta)^{-1}$ as a symmetric contravariant second order tensor acting on covariant 1-forms, we obtain
\[ (D\zeta)^{-1}((D^2\zeta)(\cdot,v,v),(D^2\zeta)(\cdot,v,v)) = 4(v^2)^T\gamma\gamma^{-1}\gamma v^2 = 4(v^2)^T\gamma v^2 = 4\gamma(v^2,v^2).
\]
Comparing this with $(D^3\zeta)(v,v,v,v) = 6\gamma(v^2,v^2)$, we finally obtain
\begin{equation} \label{zeta_rel}
(D^3\zeta)(v,v,v,v) = \frac32 (D\zeta)^{-1}((D^2\zeta)(\cdot,v,v),(D^2\zeta)(\cdot,v,v)).
\end{equation}
This relation holds at the point $e$ for arbitrary vectors $v$.

However, by Lemma \ref{lem_zeta} the form $\zeta$ and its derivatives are invariant with respect to the action of $\Pi$, and hence \eqref{zeta_rel} holds identically on ${\cal Y}$. By Lemma \ref{surf_chars_zeta} the manifold $M$ has then parallel cubic form. This completes the proof.
\end{proof}

\medskip

Let us now consider the case when $J$ is semi-simple and $\gamma = -n^{-1}g$. Then ${\cal Y}$ is an $\omega$-domain.

{\lemma \label{ze_level} Let $J$ be a real semi-simple Jordan algebra, ${\cal Y}$ the corresponding $\omega$-domain containing the unit element, $\omega$ its $\omega$-function,  and $\gamma = -n^{-1}g$ a non-degenerate trace form. Then the maximal integral manifolds of the distribution $\Delta$ from Lemma \ref{lem_zeta} are the level surfaces $\omega_c = \{ x \in {\cal Y} \,|\, \omega(x) = c \}$, $c > 0$, of the function $\omega$. }

\begin{proof}
By \eqref{der_log_det} we have $\zeta(u) = -\gamma(u,x^{-1}) = n^{-1}g(u,x^{-1}) = \frac{1}{2n}D_u(\log\det P_x)$. It follows that $\zeta$ possesses the global potential function
\begin{equation} \label{Phi_omega}
\Phi = \frac{1}{2n}\log\det P_x = \frac{1}{n}\log\omega
\end{equation}
on ${\cal Y}$. Here the second equality comes from Theorem \ref{domain_algebra}. The maximal integral manifolds of $\Delta$ are the level surfaces of $\Phi$. But these coincide with the level surfaces of $\omega$.
\end{proof}

{\theorem \label{S1} Let ${\cal Y}$ be an $\omega$-domain with function $\omega$ and pseudo-metric $\sigma$. For every $c > 0$, the level surface $\omega_c$ is a connected $(n-1)$-dimensional centro-affine submanifold of ${\cal Y}$. Its conic hull equals ${\cal Y}$, and each ray of ${\cal Y}$ intersects $\omega_c$ exactly once. As an affine hypersurface immersion, $\omega_c$ is an affine complete, Euclidean complete, homogeneous symmetric affine hypersphere with center in the origin which is asymptotic to the boundary of ${\cal Y}$ and has parallel cubic form. The restriction of the pseudo-metric $\sigma$ to $\omega_c$ is proportional to the affine metric. }

\begin{proof}
That $\omega_c$ is a centro-affine complete, homogeneous, symmetric centro-affine hypersurface, with parallel cubic form, and with conic hull equal to ${\cal Y}$, follows from Theorem \ref{alg_to_imm} and Lemma \ref{ze_level}.

That $\omega_c$ intersects each ray exactly once follows from the homogeneity of the function $\omega$.

Since $\omega = 0$ on the boundary of ${\cal Y}$ and $\omega$ is continuous on the closure of ${\cal Y}$, the level surface $\omega_c$ is closed and without boundary. It is thus Euclidean complete and asymptotic to the boundary of ${\cal Y}$.

By Lemma \ref{surf_chars_zeta} the centro-affine metric is given by the restriction of $D\zeta$ on $\omega_c$. However, $\zeta$ has the global potential $\Phi$ given by \eqref{Phi_omega}, and hence $D\zeta = \Phi'' = -n^{-1}(-\log\omega)'' = -n^{-1}\sigma$. Thus the centro-affine metric is proportional to the pseudo-metric induced by $\sigma$.

We shall now show that $\omega_c$ is an affine hypersphere. By \eqref{zeta_deriv_e} the matrix of $D\zeta$ at the point $e$ is given by $\gamma = -n^{-1}g$, and its inverse by $-ng^{-1}$. Further, $(D^2\zeta)(u,v,v) = \frac{2}{n}g(u,v^2) = \frac{2}{n}g(L_uv,v)$, and the matrix of $(D^2\zeta)(u,\cdot,\cdot)$ is given by $\frac{2}{n}L_u^Tg$. It follows that at $e$
\[ tr_{D\zeta}(D^2\zeta)(u,\cdot,\cdot) = tr(-2L_u^Tgg^{-1}) = -2tr\,L_u = -2g(u,e) = -2n\zeta(u).
\]
Now note that by Lemma \ref{lem_zeta} the form $\zeta$ and its derivatives are invariant with respect to the action of $\Pi$. Hence the above relation holds identically on ${\cal Y}$. By Lemma \ref{surf_chars_zeta} $\omega_c$ is then an affine hypersphere with center in the origin.

It follows that the Blaschke metric of $\omega_c$ is proportional to its centro-affine metric. Thus $\omega_c$ is also affine complete with respect to its Blaschke metric, and  its Blaschke metric is proportional to the pseudo-metric induced by $\sigma$.
\end{proof}

\subsection{Immersions and Jordan algebras}

Let us summarize the results of the two preceding subsections.

Let $M \subset Y$ be a non-degenerate centro-affine hypersurface with parallel cubic form, given as in Lemma \ref{surf_chars_zeta} by an integral manifold of an involutive distribution $\tilde\Delta$ on an open subset $Y \subset \mathbb R^n$. Let $\tilde\zeta$ be the closed 1-form on $Y$ satisfying $\tilde\zeta(x) \equiv 1$, where $x$ is the position vector field, such that $\tilde\Delta$ is the kernel of $\tilde\zeta$. Let $\hat D$ be the Levi-Civita connection of the pseudo-metric $D\tilde\zeta$ on $Y$, and $D$ the flat connection on $\mathbb R^n$. Choose a point $a \in M$ and consider the tensor $K = D - \hat D$ at $a$. This tensor defines a real unital Jordan algebra $J_a$ on $\mathbb R^n$ with unit element $a$. Let further $\gamma_a$ be the symmetric bilinear form defined on $J_a$ by $D\tilde\zeta$, evaluated at $a$. Then $\gamma_a$ is a non-degenerate trace form satisfying $\gamma(a,a) = -1$.

On the other hand, let $J$ be a real unital Jordan algebra of dimension $n \geq 2$, $\gamma$ a non-degenerate trace form satisfying $\gamma(e,e) = -1$, ${\cal Y}$ the connected component of the unit element $e$ in the set of invertible elements of $J$, and $\Delta$ the involutive distribution on ${\cal Y}$ from Lemma \ref{lem_zeta}. Then every maximal integral manifold of $\Delta$ is a non-degenerate centro-affine hypersurface with parallel cubic form.

We shall now consider the interplay between these relations.

{\lemma \label{imm_alg_imm} Assume above notations and set $J = J_a$, $\gamma = \gamma_a$. Then the hypersurface $M$ is an integral manifold of $\Delta$. }

\begin{proof}
Let $M'$ be the maximal integral manifold of $\Delta$ passing through the point $a$, and let $\zeta$ be the closed 1-form from Lemma \ref{lem_zeta}. Then $\Delta$ is the kernel of $\zeta$, $\tilde\Delta$ is the kernel of $\tilde\zeta$, and $M,M'$ are integral manifolds of $\tilde\Delta,\Delta$, respectively. Note that $M'$ has parallel cubic form by Theorem \ref{alg_to_imm}.

At the vector $x = a$, which is the unit element of the Jordan algebra $J_a$, we have by virtue of \eqref{zeta_def},\eqref{zeta_deriv_e} that
\begin{eqnarray*}
\zeta(u) &=& -\gamma(u,a), \\
(D\zeta)(u,u) &=& \gamma(u,u), \\
(D^2\zeta)(u,u,u) &=& -2\gamma(u,u^2)
\end{eqnarray*}
for all $u \in J_a$. On the other hand, at $a$ the values of $\tilde\zeta$ and its derivatives on $u$ are given by
\begin{eqnarray*}
\tilde\zeta(u) &=& -(D\tilde\zeta)(a,u) = -\gamma_a(a,u), \\
(D\tilde\zeta)(u,u) &=& \gamma_a(u,u), \\
(D^2\tilde\zeta)(u,u,u) &=& -2(D\tilde\zeta)(u \bullet u,u) = -2\gamma_a(u^2,u).
\end{eqnarray*}
Here we used \eqref{FderXa} in the first equality of the first line, and \eqref{K_zeta} in the first equality of the third line.

Thus at $a$ the values of $\zeta$ and its first two derivatives coincide with the values of $\tilde\zeta$ and its derivatives, respectively. It follows that the hypersurfaces $M,M'$ make a contact of order 3 at $a$. But both hypersurfaces have parallel cubic form and must hence locally coincide. Since $M'$ is centro-affine complete by Theorem \ref{alg_to_imm}, $M$ is actually contained in $M'$. The claim of the lemma now easily follows.
\end{proof}

{\theorem \label{th_main3} Let $f: M \to \mathbb R^n$ be a connected non-degenerate centro-affine hypersurface immersion with parallel cubic form. Then $f$ can be extended to a homogeneous, symmetric, affine complete, injective non-degenerate centro-affine hypersurface immersion $\bar f: \bar M \to \mathbb R^n$ with parallel cubic form. Let $Y$ be the conic hull of the image of $\bar f$. Then there exists a real unital Jordan algebra $J$ on $\mathbb R^n$ such that $Y$ is the connected component of the unit element $e$ in the set of invertible elements of $J$. Moreover, there exists a non-degenerate trace form $\gamma$ on $J$ satisfying $\gamma(e,e) = -1$ such that the immersion $\bar f$ is tangent to the distribution $\Delta$ defined on $Y$ by the pair $(J,\gamma)$ as in Lemma \ref{lem_zeta}. }

\begin{proof}
Choose an arbitrary point $\xi \in M$ and a neighbourhood $V \subset M$ of $\xi$ such that the image $f[V]$ is an embedded submanifold of $\mathbb R^n$ and each ray in $\mathbb R^n$ intersects $f[V]$ at most once. Then $U = \{ \lambda f(\varphi) \,|\, \varphi \in V,\ \lambda > 0 \}$ is an open subset of $\mathbb R^n$ which is canonically diffeomorphic to $V \times \mathbb R_{++}$. Define a function $\Phi$ on $U$ by
$\Phi(\lambda f(\varphi)) = \log\lambda$ and let $\tilde\zeta = D\Phi$.

Then $\tilde\zeta$ is a closed 1-form on $U$ and satisfies $\tilde\zeta(x) \equiv 1$, with $x$ the position vector field. Moreover, the image $f[V]$ is an integral hypersurface of the kernel $\tilde\Delta$ of $\tilde\zeta$. Put $a = f(\xi)$.

Apply Theorem \ref{main1b} to the image $f[V]$ at the point $a$, let $J$ be the Jordan algebra on $\mathbb R^n$ with unit element $a$ and $\gamma$ the trace form from this theorem. Let $\zeta$ be the 1-form \eqref{zeta_def} defined by the pair $(J,\gamma)$, and let $\Delta$ be its kernel as in Lemma \ref{lem_zeta}. Then by Lemma \ref{imm_alg_imm} the hypersurface $f[V]$ is an integral manifold of $\Delta$.

Let $M'$ be the maximal integral manifold of $\Delta$ passing through the point $a$. Then $f[V]$ is contained in $M'$. By Theorem \ref{alg_to_imm} the hypersurface $M'$ has parallel cubic form and is affine complete. Since the immersion $f$ also has parallel cubic form, and the image of $f$ shares with $M'$ the set $f[V]$, the image of $f$ must actually be contained in $M'$. The claims of the theorem now follow from the properties of $M'$ given in Theorem \ref{alg_to_imm}.
\end{proof}

Theorem \ref{th_main3} completely characterizes centro-affine hypersurface immersions with parallel cubic form. Their study is thus reduced to the study of pairs $(J,\gamma)$, where $J$ is a real unital Jordan algebra, and $\gamma$ is a non-degenerate trace form on $J$ satisfying $\gamma(e,e) = -1$. We have also the following result.

{\theorem \label{th_AHS2} Let $f: M \to \mathbb R^n$ be a proper affine hypersphere with center in the origin and with parallel cubic form. Then there exists an $\omega$-domain $Y \subset \mathbb R^n$ such that $f[M]$ is contained in some level surface $\omega_c$ of the $\omega$-function of $Y$. }

\begin{proof}
First of all, note that Theorem \ref{th_main3} is applicable. By Lemma \ref{AHS} the bilinear form $\gamma$ in Theorem \ref{th_main3} equals $-n^{-1}g$ and the algebra $J$ is semi-simple. Then by Lemma \ref{ze_level} the maximal integral manifolds of $\Delta$ are the level surfaces of the function $\omega$. This proves the theorem.
\end{proof}

Thus the study of proper affine hyperspheres with parallel cubic form can be reduced to the study of the level surfaces $\omega_c$ of $\omega$-domains, and Theorem \ref{S1} is applicable.

{\remark Lemma \ref{imm_alg_imm} essentially says that if we start with a centro-affine hypersurface immersion with parallel cubic form, pick a point on it and construct a pair $(J,\gamma)$ of a real unital Jordan algebra together with a non-degenerate trace form $\gamma$ as in Theorem \ref{main1}, and then construct an involutive distribution $\Delta$ from $(J,\gamma)$ as in Lemma \ref{lem_zeta}, whose integral manifolds have parallel cubic form, then we recover the original hypersurface. A similar result holds for the reverse way. If we start with a pair $(J,\gamma)$, construct a distribution $\Delta$ from it, choose the maximal integral manifold $M$ of $\Delta$ passing through the unit element $e$, and then construct a pair $(J',\gamma')$ from $M$ at the point $e$ as in Theorem \ref{main1}, then we recover the original pair $(J,\gamma)$. }

Consider a centro-affine hypersurface immersion $f: M \to \mathbb R^n$ with parallel cubic form. In order to construct a pair $(J,\gamma)$ that in turn yields the centro-affine completion of the immersion $f$, we had to choose a point $y \in {\cal D} = M \times \mathbb R_{++}$. The next result shows that the Jordan algebra $J$ is essentially independent of the basepoint $y$.

{\lemma Let $f: M \to \mathbb R^n$ be a non-degenerate centro-affine hypersurface immersion with parallel cubic form. Let $y,y' \in {\cal D} = M \times \mathbb R_{++}$ be points and let $J_y,J_{y'}$ be the Jordan algebras defined by the difference tensor $K = D - \hat D$ on the tangent spaces $T_y{\cal D},T_{y'}{\cal D}$, respectively. Then $J_y,J_{y'}$ are isomorphic. }

\begin{proof}
Choose a smooth path in ${\cal D}$ connecting $y,y'$ and transport the tangent space $T_y{\cal D}$ along this path using the parallel transport of the connection $\hat D$. In this way we obtain a map $A_{y,y'}: T_y{\cal D} \to T_{y'}{\cal D}$, which by Lemma \ref{K_parallel} preserves the tensor $K$. Hence the map $A_{y,y'}$ is an isomorphism between $J_y$ and $J_{y'}$.
\end{proof}

The isomorphism $A_{y,y'}$ may not be canonical, however, as it in general depends on the path linking $y$ and $y'$. In particular, any closed path leading back to the original basepoint $y$ induces an automorphism of $J_y$. Clearly this automorphism has to preserve the bilinear form $\gamma$ determined by the immersion.

In any case the pairs $(J,\gamma)$ constructed from different points on the immersion lead to the same connected component ${\cal Y}$, because the centro-affine completion of the immersion is unique and ${\cal Y}$ is the conic hull of its image.

\section{Classification} \label{sec_classification}

In this section we give a complete classification of the non-degenerate proper affine hyperspheres with parallel cubic form. It is based on the classification of finite-dimensional simple real Jordan algebras and the fact that a semi-simple Jordan algebra is the direct product of uniquely determined simple factors. We also classify those centro-affine hypersurface immersions with parallel cubic form whose associated Jordan algebra is semi-simple.

\subsection{Decomposition}

In this subsection we investigate which impact the decomposition of a unital Jordan algebra $J$ into a direct sum of ideals $J_1,\dots,J_r$ has on the centro-affine hypersurfaces with parallel cubic form which $J$ defines.

For an element $x \in J$, we shall denote by $x_k$ its projection on the ideal $J_k$. Thus $x = \sum_{k=1}^r x_k$ with $x_k \in J_k$. The projections $e_k$ of the unit element $e$ are the unit elements in the ideals $J_k$, and all $J_k$ are also unital. Pass to a coordinate system which is adapted to the decomposition $J = \oplus_{k=1}^r J_k$. Since the $J_k$ are ideals, all operators $L_x$ are block-diagonal. Namely, we have $L_x = \diag(L_{x_1},\dots,L_{x_r})$, where $L_{x_k}$ is the operator of multiplication with $x_k$ in $J_k$. This implies that $(x^2)_k = x_k^2$, and the operators $P_x$ are also block-diagonal. Namely, $P_x = \diag(P_{x_1},\dots,P_{x_r})$, where $P_{x_k}$ is the quadratic operator of $x_k$ in $J_k$. It follows that $x$ is invertible in $J$ if and only if all $x_k$ are invertible in $J_k$, and the inverse of $x$ is given by $\sum_{k=1}^r x_k^{-1}$ (of course, $x_k$ is not invertible in $J$ if $r > 1$, so the inverse of $x_k$ is always assumed in $J_k$). Let ${\cal Y}$ be the connected component of $e$ in the set of invertible elements in $J$, and ${\cal Y}_k$ the connected component of $e_k$ in the set of invertible elements of $J_k$. From the above it then follows that ${\cal Y}$ is the sum of the ${\cal Y}_k$, ${\cal Y} = \{ y \,|\, y_k \in {\cal Y}_k\ \forall\ k = 1,\dots,r \}$.

{\lemma \label{decomp_nondeg} Assume above notations. Let $\gamma$ be a trace form on $J$, and denote by $\gamma_k$ the restriction of $\gamma$ to the subspace $J_k$. Then $\gamma(x,y) = \sum_{k=1}^r \gamma_k(x_k,y_k)$ and $\gamma$ is non-degenerate if and only if all $\gamma_k$ are non-degenerate. Moreover, $\gamma_k$ is a trace form on $J_k$ for all $k$. On the other hand, if $\sigma_k$ are trace forms on $J_k$, then the form $\sigma$ on $J$ given by $\sigma(x,y) = \sum_{k=1}^r \sigma_k(x_k,y_k)$ is a trace form. }

\begin{proof}
For $k \not= l$ we have $\gamma(x_k,y_l) = \gamma(x_k \bullet e,y_l) = \gamma(e,x_k \bullet y_l) = \gamma(e,0) = 0$. Thus
\[ \gamma(x,y) = \sum_{k,l=1}^r \gamma(x_k,y_l) = \sum_{k=1}^r \gamma(x_k,y_k) = \sum_{k=1}^r \gamma_k(x_k,y_k).
\]
The matrix of $\gamma$ is then given by $\diag(\gamma_1,\dots,\gamma_r)$, and it is invertible if and only if all blocks $\gamma_k$ are invertible. A symmetric block-diagonal matrix $\sigma = \diag(\sigma_1,\dots,\sigma_r)$ satisfies $L_x^T\sigma = \sigma L_x$ for all $x \in J$, i.e., $\sigma$ is a trace form on $J$, if and only if $L_{x_k}^T\sigma_k = \sigma_k L_{x_k}$ for all $k$ and all $x_k \in J_k$, i.e., if all $\sigma_k$ are trace forms on $J_k$.
\end{proof}

{\lemma \label{sum_zeta_Phi} Assume the notations at the beginning of the subsection. Let $\gamma$ be a non-degenerate trace form on $J$ and $\gamma_k$ its restrictions to $J_k$. Let $\zeta,\zeta_k$ be the 1-forms \eqref{zeta_def} defined on ${\cal Y},{\cal Y}_k$ by the pairs $(J,\gamma),(J_k,\gamma_k)$, respectively. Then we have $\zeta(u) = \sum_{k=1}^r \zeta_k(u_k)$ for all $u \in J$. Let $x \in {\cal Y}$ and let $U$ be a neighbourhood of $x$ such that there exists a local potential $\Phi: U \to \mathbb R$ of $\zeta$. Then there exist neighbourhoods $U_k$ of $x_k$ in ${\cal Y}_k$ and local potentials $\Phi_k: U_k \to \mathbb R$ of $\zeta_k$ such that $U' = U_1 + \dots + U_r \subset U$ and $\Phi(y) = \sum_{k=1}^r \Phi_k(y_k)$ for all $y \in U'$. }

\begin{proof}
We have
\[ \zeta(u) = -\gamma(u,x^{-1}) = -\sum_{k=1}^r \gamma_k(u_k,x_k^{-1}) = \sum_{k=1}^r \zeta_k(u_k).
\]
Here the second equality comes from Lemma \ref{decomp_nondeg}. Note that the trace forms $\gamma_k$ are non-degenerate by Lemma \ref{decomp_nondeg}, and $\zeta,\zeta_k$ are closed by Lemma \ref{lem_zeta}. Let us choose connected neighbourhoods $U_k$ of $x_k$ such that there exist local potentials $\Phi_k: U_k \to \mathbb R$ of $\zeta_k$ and $U' = U_1 + \dots + U_r \subset U$. By possibly adding a constant to one of the $\Phi_k$ we may assume that $\sum_{k=1}^r \Phi_k(x_k) = \Phi(x)$. Let $y \in U'$ be an arbitrary point. Choose smooth paths $\sigma_k: [0,1] \to U_k$ connecting $x_k$ with $y_k$. Then the path $\sigma$ defined by $\sigma(\tau) = \sum_{k=1}^r \sigma_k(\tau)$ lies in $U$ and connects $x$ with $y$. We have
\[ \Phi(y) = \Phi(x) + \int_0^1 \zeta(\dot\sigma(\tau)) d\tau = \sum_{k=1}^r \Phi_k(x_k) + \int_0^1 \sum_{k=1}^r \zeta_k(\dot\sigma_k(\tau)) d\tau = \sum_{k=1}^r \Phi_k(y_k).
\]
This completes the proof.
\end{proof}

By Theorem \ref{th_main3} centro-affine hypersurface immersions with parallel cubic form can be characterized as level surfaces of a potential $\Phi$ of the form $\zeta$ defined as in \eqref{zeta_def} by some real unital Jordan algebra $J$ and a non-degenerate trace form $\gamma$ satisfying $\gamma(e,e) = -1$. Lemma \ref{sum_zeta_Phi} then implies that in order to describe such an immersion with decomposable Jordan algebra $J$, one only needs to compute the potentials $\Phi_k$ for the indecomposable factors $J_k$ of $J$. Note that if $\gamma$ is a non-degenerate trace form on $J$ satisfying $\gamma(e,e) = -1$, then its restrictions $\gamma_k$ are non-degenerate trace forms on $J_k$, but they do not need in general to satisfy the condition $\gamma_k(e_k,e_k) \not= 0$. If $\gamma_k(e_k,e_k) = 0$, however, then Lemma \ref{Gamma_z} is no more valid and the maximal integral manifolds of the kernel $\Delta_k$ of $\zeta_k$ are no more centro-affine. This is why a centro-affine hypersurface immersion with parallel cubic form which corresponds to a decomposable Jordan algebra does not need be itself decomposable.

{\lemma \label{reverse_ca} Assume the notations at the beginning of the subsection. Let $\gamma_k$ be non-degenerate trace forms on $J_k$ such that the numbers $\gamma_k(e_k,e_k)$ are not all zero, and let $\zeta_k$ be the closed 1-forms defined as in \eqref{zeta_def} by $\gamma_k$ on ${\cal Y}_k$. Let $c_k$ be nonzero real numbers such that $\sum_{k=1}^r c_k \gamma_k(e_k,e_k) \not= 0$. Let $U_k \subset {\cal Y}_k$ be neighbourhoods such that there exist potentials $\Phi_k: U_k \to \mathbb R$ of $\zeta_k$. Define a function $\Phi(x) = \sum_{k=1}^r c_k \Phi_k(x_k)$ on $U_1 + \dots + U_r \subset {\cal Y}$.

Then the level surfaces of $\Phi$ are centro-affine hypersurfaces with parallel cubic form. }

\begin{proof}
By possibly multiplying all $c_k$ and $\Phi$ with the same nonzero constant, we can assume without loss of generality that $\sum_{k=1}^r c_k \gamma_k(e_k,e_k) = -1$.

Define a symmetric bilinear form $\gamma$ on $J$ by $\gamma(x,y) = \sum_{k=1}^r c_k \gamma_k(x_k,y_k)$. Then by the last part of Lemma \ref{decomp_nondeg} $\gamma$ is a non-degenerate trace form, and it satisfies $\gamma(e,e) = -1$.

Let $\zeta$ be the 1-form defined on ${\cal Y}$ by $(J,\gamma)$ as in \eqref{zeta_def}. From Lemma \ref{sum_zeta_Phi} it follows that $\zeta(u) = \sum_{k=1}^r c_k \zeta_k(u_k)$. Then the relations $\zeta_k = D\Phi_k$ imply that $\zeta = D\Phi$, and $\Phi$ is a potential of $\zeta$. The claim of the lemma now follows from Theorem \ref{alg_to_imm}.
\end{proof}

Analogous to the fifth relation in Lemma \ref{FderX} it can be proven that the logarithmic homogeneity constant of $\Phi$ is given by $-\gamma(e,e)$. Hence the condition $\sum_{k=1}^r c_k \gamma_k(e_k,e_k) \not= 0$ is equivalent to the condition that this constant is nonzero.

\medskip

Let us now consider the case when the centro-affine hypersurface with parallel cubic form is an affine hypersphere with center in the origin. By Lemma \ref{AHS} this situation occurs if and only if $J$ is semi-simple and $\gamma = -n^{-1}g$, where $g$ is given by \eqref{bilinear}. By Theorem \ref{semi_simple_decomp} $J$ decomposes into a direct sum of simple ideals. We shall show that this decomposition induces a representation of the affine hypersphere as a Calabi product of lower-dimensional affine hyperspheres with parallel cubic form. Note that in Subsection \ref{sec_imm_alg} we excluded Jordan algebras of dimension 1. The simple factors of a semi-simple algebra $J$ may, however, have dimension 1. Therefore in what follows below we consider the point as an affine hypersphere with parallel cubic form and include it in the list of possible factors occurring in the Calabi product.

{\lemma \label{lem_Calabi} Let $J = \oplus_{k=1}^r J_k$ be the decomposition of a real semi-simple Jordan algebra into simple factors. Let ${\cal Y},{\cal Y}_k$ be the $\omega$-domains of $J,J_k$ containing the respective unit element. Let $c > 0,c_k > 0$, $k = 1,\dots,r$, be constants. Let $\omega_c,\omega_{k,c_k}$ be the corresponding level surfaces of the $\omega$-functions $\omega,\omega_k$ on ${\cal Y},{\cal Y}_k$, respectively. Then $\omega_c$ is a Calabi product of the level surfaces $\omega_{k,c_k}$, which may include points as factors. }

\begin{proof}
Let $n_k$ be the dimension of $J_k$, and $n = \sum_{k=1}^r n_k$ the dimension of $J$. By Theorem \ref{S1} all level surfaces $\omega_{k,c_k}$ with $n_k \geq 2$ are indeed proper affine hyperspheres with center in the origin. The level surfaces $\omega_{k,c_k}$ for which $n_k = 1$ are points. Let $\omega_{k,c_k}$ be given as an affine hypersurface immersion by the inclusion map $i_k: \omega_{k,c_k} \to J_k$. Define the $(r-1)$-dimensional affine subspace $A \subset \mathbb R^r$ by
\[ A = \left\{ t = (t_1,\dots,t_r)^T \,|\, \sum_{k=1}^r n_kt_k = \log c - \sum_{k=1}^r \log c_k \right\}.
\]
Let us define the immersion $f: A \times \prod_{k=1}^r \omega_{k,c_k} \to J$ by
\begin{equation} \label{Calabi_product}
f(t,x_1,\dots,x_r) = \sum_{k=1}^r e^{t_k} i_k(x_k).
\end{equation}
This is a Calabi product of the affine hyperspheres $\omega_{k,c_k}$ \cite{XLi93}.

Let us show that its image is the level surface $\omega_c$. Recall that the functions $\omega,\omega_k$ are homogeneous of degree $n,n_k$, respectively. Moreover, $\det P_x = \prod_{k=1}^r \det P_{x_k}$ and hence by Theorem \ref{domain_algebra} $\omega(x) = \prod_{k=1}^n \omega_k(x_k)$ for all $x \in {\cal Y}$. We then have
\begin{eqnarray*}
\omega(f(t,x_1,\dots,x_r)) &=& \omega\left( \sum_{k=1}^r e^{t_k} x_k \right) = \prod_{k=1}^r \omega_k\left( e^{t_k} x_k \right) = \prod_{k=1}^r e^{n_kt_k} \omega_k(x_k) = e^{\sum_{k=1}^r n_kt_k} \prod_{k=1}^r c_k \\ &=& e^{\log c - \sum_{k=1}^r \log c_k} \prod_{k=1}^r c_k = c.
\end{eqnarray*}
Hence the image of the immersion $f$ is contained in the level surface $\omega_c$.

On the other hand, let $x = \sum_{k=1}^r x_k \in \omega_c$, with $x_k \in {\cal Y}_k$. Since the conic hull of the level surface $\omega_{k,c_k}$ is ${\cal Y}_k$, there exist $\alpha_k > 0$ and $y_k \in \omega_{k,c_k}$ such that $x_k = \alpha_k y_k$. We then have $\omega_k(x_k) = \alpha_k^{n_k} \omega_k(y_k) = \alpha_k^{n_k} c_k$. Since $\omega(x) = c$, we obtain
\[ \log c = \log\omega(x) = \log\prod_{k=1}^r \omega_k(x_k) = \sum_{k=1}^r \log\left( \alpha_k^{n_k} c_k \right) = \sum_{k=1}^r \left( n_k\log\alpha_k + \log c_k \right).
\]
Hence $t = (\log\alpha_1,\dots,\log\alpha_r)^T \in A$, and $f(t,y_1,\dots,y_r) = \sum_{k=1}^r \alpha_k y_k = x$, which proves that $\omega_c$ is contained in the image of $f$.
\end{proof}

{\remark Originally, the Calabi product of proper affine hyperspheres was defined for two factors. It can, however, in a straightforward manner be extended to multiple factors, which was accomplished in \cite{XLi93} for the hyperbolic case. The formula \eqref{Calabi_product} is equivalent to, but much simpler than the existing definitions, and is valid for affine hyperspheres of arbitrary signature. }

Thus, in order to classify the proper affine hyperspheres with parallel cubic form, it is sufficient to classify those which are defined by the $\omega$-domains of simple Jordan algebras.

\subsection{Real and complex Jordan algebras}

In this subsection we consider the situation when a real semi-simple Jordan algebra $J_{\mathbb R}$ of dimension $2n$ is isomorphic to a complex Jordan algebra $J_{\mathbb C}$ of dimension $n$. This is motivated by Lemma \ref{lem_central_simple}, which implies that this case occurs whenever a simple real Jordan algebra is not central-simple.

For convenience, we identify the underlying vector spaces such that the multiplication $\bullet$ in $J_{\mathbb R}$ and $J_{\mathbb C}$ is the same operation. It follows that the linear mapping $P_x: y \mapsto 2x \bullet (x \bullet y) - (x \bullet x) \bullet y$ is also the same operation in $J_{\mathbb R}$ and $J_{\mathbb C}$. If $x_{\mathbb C} \in \mathbb C^n$ is the coordinate vector of some point $x$ in $J_{\mathbb C}$, then we assign to $x$ the coordinate vector $x_{\mathbb R} = \begin{pmatrix} Re\,x \\ Im\,x \end{pmatrix} \in \mathbb R^{2n}$ in $J_{\mathbb R}$. The linear operators $L_x,P_x$ can then be represented as complex matrices $L_x^{\mathbb C},P_x^{\mathbb C} \in \mathbb C^{n \times n}$ when acting on complex coordinate vectors, and as real matrices $L_x^{\mathbb R},P_x^{\mathbb R} \in \mathbb R^{2n \times 2n}$ when acting on real coordinate vectors.

{\lemma \label{L_relation} The matrices $L_x^{\mathbb R},P_x^{\mathbb R}$ are given by
\[ L_x^{\mathbb R} = \begin{pmatrix} Re\,L_x^{\mathbb C} & -Im\,L_x^{\mathbb C} \\ Im\,L_x^{\mathbb C} & Re\,L_x^{\mathbb C} \end{pmatrix}, \qquad P_x^{\mathbb R} = \begin{pmatrix} Re\,P_x^{\mathbb C} & -Im\,P_x^{\mathbb C} \\ Im\,P_x^{\mathbb C} & Re\,P_x^{\mathbb C} \end{pmatrix}.
\] }

\begin{proof}
Let $y$ be an arbitrary element. We then have
\[ (L_xy)_{\mathbb C} = L_x^{\mathbb C}y_{\mathbb C} = (Re\,L_x^{\mathbb C}Re\,y_{\mathbb C} - Im\,L_x^{\mathbb C}Im\,y_{\mathbb C}) + i(Re\,L_x^{\mathbb C}Im\,y_{\mathbb C} + Im\,L_x^{\mathbb C}Re\,y_{\mathbb C}).
\]
Hence
\[ (L_xy)_{\mathbb R} = \begin{pmatrix} Re\,L_x^{\mathbb C}Re\,y_{\mathbb C} - Im\,L_x^{\mathbb C}Im\,y_{\mathbb C} \\ Re\,L_x^{\mathbb C}Im\,y_{\mathbb C} + Im\,L_x^{\mathbb C}Re\,y_{\mathbb C} \end{pmatrix} = \begin{pmatrix} Re\,L_x^{\mathbb C} & -Im\,L_x^{\mathbb C} \\ Im\,L_x^{\mathbb C} & Re\,L_x^{\mathbb C} \end{pmatrix}\begin{pmatrix} Re\,y_{\mathbb C} \\ Im\,y_{\mathbb C} \end{pmatrix} = L_x^{\mathbb R}y_{\mathbb R}.
\]
The proof for the matrix $P_x^{\mathbb R}$ is similar.
\end{proof}

Note that $P_x^{\mathbb R} = \begin{pmatrix} Re\,P_x^{\mathbb C} & -Im\,P_x^{\mathbb C} \\ Im\,P_x^{\mathbb C} & Re\,P_x^{\mathbb C} \end{pmatrix}$ can be written as $\frac12 \begin{pmatrix} iI & I \\ I & iI \end{pmatrix} \begin{pmatrix} P_x^{\mathbb C} & 0 \\ 0 & \bar P_x^{\mathbb C} \end{pmatrix} \begin{pmatrix} -iI & I \\ I & -iI \end{pmatrix}$. We then have $\det P_x^{\mathbb R} = |\det P_x^{\mathbb C}|^2$, and the $\omega$-function on the $\omega$-domains of $J_{\mathbb R}$ is given by
\begin{equation} \label{omega_complex}
\omega(x) = |\det P_x^{\mathbb C}|.
\end{equation}
From Lemma \ref{L_relation} we also get the following result.

{\corollary \label{g_t_corr} The linear forms $t_{\mathbb C},t_{\mathbb R}$ given by \eqref{trace_form} in the Jordan algebras $J_{\mathbb C},J_{\mathbb R}$, respectively, are related by $t_{\mathbb R}(u) = 2Re\,t_{\mathbb C}(u)$. The bilinear symmetric forms $g_{\mathbb C},g_{\mathbb R}$ given by \eqref{bilinear} in the Jordan algebras $J_{\mathbb C},J_{\mathbb R}$, respectively, are related by $g_{\mathbb R}(u,v) = 2Re\,g_{\mathbb C}(u,v)$. }

\begin{proof}
We have
\[ t_{\mathbb R}(u) = tr\,L_u^{\mathbb R} = 2tr\,Re\,L_u^{\mathbb C} = 2Re\,tr\,L_u^{\mathbb C} = 2Re\,t_{\mathbb C}(u).
\]
It follows that
\[ g_{\mathbb R}(u,v) = t_{\mathbb R}(u \bullet v) = 2Re\,t_{\mathbb C}(u \bullet v) = 2Re\,g_{\mathbb C}(u,v). \qedhere
\]
\end{proof}

The matrix of $g_{\mathbb R}$ can then be written as
\[ 2\begin{pmatrix} Re\,g_{\mathbb C} & -Im\,g_{\mathbb C} \\ -Im\,g_{\mathbb C} & -Re\,g_{\mathbb C} \end{pmatrix} = \begin{pmatrix} iI & I \\ -I & -iI \end{pmatrix} \begin{pmatrix} g_{\mathbb C} & 0 \\ 0 & \bar g_{\mathbb C} \end{pmatrix} \begin{pmatrix} -iI & I \\ I & -iI \end{pmatrix}.
\]
It follows that $\det g_{\mathbb R} = (-4)^n|\det g_{\mathbb C}|^2$, and $J_{\mathbb R}$ is semi-simple if and only if $J_{\mathbb C}$ is.

{\lemma Let $c \in \mathbb C$ be a nonzero number. Then the $\mathbb R$-bilinear form $\gamma(u,v) = 2Re(cg_{\mathbb C}(u,v))$ is a non-degenerate trace form on $J_{\mathbb R}$. }

\begin{proof}
First note that $\gamma$ is symmetric.

By the $\mathbb C$-bilinearity of $g_{\mathbb C}$ we have $\gamma(u,v) = 2Re\,g_{\mathbb C}(cu,v) = g_{\mathbb R}(cu,v)$, where the second equality comes from the preceding corollary. Suppose that there exists a nonzero $u$ such that $\gamma(u,v) = 0$ for all $v$. Then $cu$ is also nonzero and hence $g_{\mathbb R}$ is degenerate, which contradicts the semi-simplicity of $J_{\mathbb R}$. Hence $\gamma$ is non-degenerate.

Finally,
\[ \gamma(u \bullet v,w) = 2Re(cg_{\mathbb C}(u \bullet v,w)) = 2Re(cg_{\mathbb C}(u,v \bullet w)) = \gamma(u,v \bullet w),
\]
where the second equality follows from \eqref{3symmetry}. Hence $\gamma$ is a trace form.
\end{proof}

The 1-form \eqref{zeta_def} for this particular choice of the trace form $\gamma$ is given by $\zeta(u) = -2Re(cg_{\mathbb C}(u,x^{-1}))$ at the point $x$. By \eqref{der_log_det} we have $D_u(\log\det P_x^{\mathbb C}) = 2g_{\mathbb C}(u,x^{-1})$, and hence
\[ \zeta = -Re\left(cD(\log\det P_x^{\mathbb C})\right) = -D\left( Re\left(c\cdot\log\det P_x^{\mathbb C}\right) \right).
\]
The form $\zeta$ has therefore a potential given by
\begin{equation} \label{complex_potential}
\Phi(x) = -Re\left(c\cdot\log\det P_x^{\mathbb C}\right).
\end{equation}
Note that since the complex logarithm is multi-valued, this potential can in general only be defined locally.

\subsection{Complex Jordan algebras}

In this subsection we provide the classification of all complex simple Jordan algebras $J$ and compute their determinant $\det P_x^{\mathbb C}$. Recall that an element $x$ is invertible if and only if $\det P_x^{\mathbb C} = 0$. Hence the set of singular elements is a variety of real codimension 2, and the connection component ${\cal Y}$ of the unit element in the set of invertible elements is open and dense in $J$. By Lemma \ref{lem_central_simple} $J$ is also central-simple, and the invertible central elements $z$ of $J$ have the form $ce$ for some nonzero $c \in \mathbb C$. If we now consider $J$ as an algebra over $\mathbb R$, then its bilinear form \eqref{bilinear} is a non-degenerate trace form, and by Theorem \ref{form_exist} any non-degenerate trace form $\gamma$ on $J$ must have the form $\gamma(u,v) = g_{\mathbb R}(z \bullet u,v) = g_{\mathbb R}(cu,v) = 2Re(cg_{\mathbb C}(u,v))$ for some nonzero $c \in \mathbb C$, where the last equality comes from Corollary \ref{g_t_corr}. The 1-form $\zeta$ defined by the pair $(J,\gamma)$ then locally has a potential of the form \eqref{complex_potential}.

{\theorem \cite[p.66--68]{McCrimmon} Let $J$ be a finite-dimensional simple Jordan algebra over $\mathbb C$. Then $J$ is exactly one of the following:

\begin{itemize}
\item $\mathbb C$,
\item ${\cal J}ord_m(I)$, the $m$-dimensional complex quadratic factor for $m \geq 3$,
\item $S_m(\mathbb C)$, the algebra of complex symmetric $m \times m$ matrices for $m \geq 3$,
\item $M_m(\mathbb C)$, the algebra of all complex $m \times m$ matrices for $m \geq 3$,
\item $H_m(Q,\mathbb C)$, the algebra of $m \times m$ Hermitian matrices with entries being split quaternions over $\mathbb C$ for $m \geq 3$,
\item $H_3(O,\mathbb C)$, the algebra of $3 \times 3$ Hermitian matrices with entries being split octonions over $\mathbb C$.
\end{itemize} }

\subsubsection{Ground field $\mathbb C$}

Since $\mathbb C$ has a commutative multiplication, it is in particular a Jordan algebra. Its unit element is 1, and all nonzero elements are invertible. For $x \in \mathbb C$, we have $L_x^{\mathbb C} = x$, $P_x^{\mathbb C} = x^2$. Hence $\det P_x^{\mathbb C} = x^2$, and by \eqref{omega_complex} $\omega(x) = |x|^2$.

\subsubsection{Quadratic factor ${\cal J}ord_m(I)$}

We shall describe the quadratic factor for a general quadratic form for future reference. Let $Q$ be a non-degenerate symmetric bilinear form on $\mathbb C^m$ and $e \in \mathbb C^m$ a distinguished point such that $e^TQe = 1$. Then the quadratic factor ${\cal J}ord_m(Q,e)$ is the vector space $\mathbb C^m$ equipped with the product \cite[p.75]{McCrimmon}
\[ x \bullet y = e^TQx \cdot y + e^TQy \cdot x - x^TQy \cdot e.
\]
The operators of multiplication are hence given by $L_x^{\mathbb C} = e^TQx \cdot I + xe^TQ - ex^TQ$. A straightforward computation yields $P_x^{\mathbb C} = x^TQx \cdot I + 4e^TQx \cdot xe^TQ - 2xx^TQ - 2x^TQx \cdot ee^TQ$ and $\det P_x^{\mathbb C} = (x^TQx)^m$. A point $x$ is invertible if and only if $x^TQx \not= 0$.

Note that over $\mathbb C$ every non-degenerate symmetric bilinear form is equivalent to the identity matrix, and we can set $Q = I$ without loss of generality. We then have
\[ \det P_x^{\mathbb C} = (x^Tx)^m,\qquad \omega(x) = |x^Tx|^m.
\]

\subsubsection{Symmetric matrices $S_m(\mathbb C)$}

The Jordan algebra $S_m(\mathbb C)$ is the space of complex symmetric $m \times m$ matrices equipped with the multiplication $A \bullet B = \frac{AB+BA}{2}$ \cite[pp.58--59; pp.66--68]{McCrimmon}. Then $P_AB = ABA$ \cite[p.82]{McCrimmon}, and $A$ is invertible in the Jordan algebra if and only if it is invertible as a matrix. Let us compute the determinant $\det P_A^{\mathbb C}$. Assume that the matrix $A$ is diagonalizable, with eigenvalues $\lambda_1,\dots,\lambda_m$ and corresponding eigenvectors $v_1,\dots,v_m$. Then the matrices $B_{ij} = v_iv_j^T + v_jv_i^T$, $1 \leq i \leq j \leq m$, form a basis of $S_m(\mathbb C)$. We have
\[ P_AB_{ij} = A(v_iv_j^T + v_jv_i^T)A = \lambda_i\lambda_j(v_iv_j^T + v_jv_i^T) = \lambda_i\lambda_jB_{ij}.
\]
Hence the $B_{ij}$ are eigenvectors of the operator $P_A$ with eigenvalues $\lambda_i\lambda_j$. Therefore we obtain
\[ \det P_A^{\mathbb C} = \prod_{i \leq j} \lambda_i\lambda_j = \prod_i \lambda_i^{m+1} = (\det A)^{m+1}.
\]
Now note that the set of diagonalizable matrices is open and dense in $S_m(\mathbb C)$. Since $\det P_A^{\mathbb C}$ is a continuous function of $A$, the above formula must be valid for all matrices. We then obtain $\omega(A) = |\det A|^{m+1}$.

\subsubsection{Full matrices $M_m(\mathbb C)$}

The Jordan algebra $M_m(\mathbb C)$ is the space of complex $m \times m$ matrices equipped with the multiplication $A \bullet B = \frac{AB+BA}{2}$ \cite[pp.58--59]{McCrimmon}. Again $P_AB = ABA$, and $A$ is invertible in the Jordan algebra if and only if it is invertible as a matrix. Assume that the matrix $A$ is diagonalizable, with eigenvalues $\lambda_1,\dots,\lambda_m$. Let the corresponding eigenvectors of $A$ be $v_1,\dots,v_m$, and those of $A^T$ be $w_1,\dots,w_m$. Then the matrices $B_{ij} = v_iw_j^T$, $1 \leq i,j \leq m$, form a basis of $M_m(\mathbb C)$. We have
\[ P_AB_{ij} = Av_iw_j^TA = \lambda_i\lambda_jv_iw_j^T = \lambda_i\lambda_jB_{ij}.
\]
As above we obtain
\[ \det P_A^{\mathbb C} = \prod_{i,j} \lambda_i\lambda_j = \prod_i \lambda_i^{2m} = (\det A)^{2m}.
\]
As above, this is valid for all matrices $A$. We then obtain $\omega(A) = |\det A|^{2m}$.

\subsubsection{Split quaternionic matrices $H_m(Q,\mathbb C)$}

The split quaternions over $\mathbb C$ can be represented by complex $2 \times 2$ matrices, with conjugation given by \cite[p.66]{McCrimmon}
\begin{equation} \label{split_quat_conj}
conj\,\begin{pmatrix} a & b \\ c & d \end{pmatrix} = \begin{pmatrix} d & -b \\ -c & a \end{pmatrix}.
\end{equation}
Hence the space $H_m(Q,\mathbb C)$ of Hermitian matrices with split quaternionic entries can be represented as the space of $2m \times 2m$ complex skew-Hamiltonian matrices $\begin{pmatrix} A & B \\ C & A^T \end{pmatrix}$, where $B,C$ are skew-symmetric. If we introduce the matrix $J = \begin{pmatrix} 0 & I \\ -I & 0 \end{pmatrix}$, then we can represent the skew-Hamiltonian matrices by products $JS$, where $S$ is a skew-symmetric $2m \times 2m$ matrix \cite[Proposition 1]{Ikramov01}.

The Jordan algebra $H_m(Q,\mathbb C)$ can be represented as the space of complex skew-Hamiltonian $2m \times 2m$ matrices equipped with the multiplication $A \bullet B = \frac{AB+BA}{2}$ (cf.\ \cite[p.58]{McCrimmon}). Substituting $A = JS$, $B = JT$, $A \bullet B = J(S \bullet T)$, we can view it equivalently as the space $A_{2m}(\mathbb C)$ of complex skew-symmetric $2m \times 2m$ matrices equipped with the multiplication
\begin{equation} \label{split_quat_mult}
S \bullet T = \frac{SJT+TJS}{2}.
\end{equation}
We shall adopt this latter point of view. Then it is straightforward to show that the operator $P_S^{\mathbb C}$ acts like $T \mapsto SJTJS$. Let us compute its determinant. Introduce the operator ${\bf J}$ acting on $A_{2m}(\mathbb C)$ like $T \mapsto JTJ$. It is not hard to see that this operator has two eigenvalues $+1,-1$ with multiplicities $m(m-1),m^2$, respectively. Its determinant hence equals $(-1)^{m^2}$. Suppose that $S$ is diagonalizable. Let $\lambda_1,\dots,\lambda_{2m}$ be its eigenvalues and $v_1,\dots,v_{2m}$ the corresponding eigenvectors. Then the matrices $T_{ij} = v_iv_j^T - v_jv_i^T$, $1 \leq i < j \leq 2m$, constitute a basis of $A_{2m}(\mathbb C)$. We have
\[ P_S^{\mathbb C} \circ {\bf J}: T_{ij} \mapsto SJ^2(v_iv_j^T - v_jv_i^T)J^2S = S(v_iv_j^T - v_jv_i^T)S = -\lambda_i\lambda_j(v_iv_j^T - v_jv_i^T) = -\lambda_i\lambda_jT_{ij},
\]
and $T_{ij}$ is an eigenvector of $P_S^{\mathbb C} \circ {\bf J}$ with eigenvalue $-\lambda_i\lambda_j$. We then get
\[ \det P_S^{\mathbb C} \cdot \det {\bf J} = \det (P_S^{\mathbb C} \circ {\bf J}) = \prod_{i < j}(-\lambda_i\lambda_j) = (-1)^{m(2m-1)}\prod_i \lambda_i^{2m-1} = (-1)^{m(2m-1)}(\det S)^{2m-1}.
\]
This yields $\det P_S^{\mathbb C} = (-1)^{m(2m-1)-m^2}(\det S)^{2m-1} = (\det S)^{2m-1} = (\pf S)^{2(2m-1)}$, where $\pf S$ is the Pfaffian of $S$. As above, this must hold for all $S$. We also see that $S$ is invertible in the Jordan algebra if and only if it is invertible as a matrix. We then obtain $\omega(S) = |\pf S|^{2(2m-1)}$.

\subsubsection{Split octonionic matrices $H_3(O,\mathbb C)$} \label{subsec_split_oct}

The complex split octonions $O$ are an 8-dimensional non-commutative, non-associative algebra over $\mathbb C$ which is generated by three hypercomplex units $j,k,l$ with multiplication table \cite[pp.64--66]{McCrimmon}

\smallskip

\begin{tabular}{|c|c|c|c|c|c|c|c|}
\hline
1 & $j$ & $k$ & $jk$ & $l$ & $jl$ & $kl$ & $(jk)l$ \\
\hline
$j$ & 1 & $jk$ & $k$ & $jl$ & $l$ & $-(jk)l$ & $-kl$ \\
\hline
$k$ & $-jk$ & 1 & $-j$ & $kl$ & $(jk)l$ & $l$ & $jl$ \\
\hline
$jk$ & $-k$ & $j$ & $-1$ & $(jk)l$ & $kl$ & $-jl$ & $-l$ \\
\hline
$l$ & $-jl$ & $-kl$ & $-(jk)l$ & 1 & $-j$ & $-k$ & $-jk$ \\
\hline
$jl$ & $-l$ & $-(jk)l$ & $-kl$ & $j$ & $-1$ & $jk$ & $k$ \\
\hline
$kl$ & $(jk)l$ & $-l$ & $jl$ & $k$ & $-jk$ & $-1$ & $-j$ \\
\hline
$(jk)l$ & $kl$ & $-jl$ & $l$ & $jk$ & $-k$ & $j$ & 1 \\
\hline
\end{tabular}

\smallskip

The {\sl conjugate} of $a = c_1+c_2j+c_3k+c_4jk+c_5l+c_6jl+c_7kl+c_8(jk)l$ is given by $\bar a = c_1-c_2j-c_3k-c_4jk-c_5l-c_6jl-c_7kl-c_8(jk)l$ and its {\sl norm} is the complex number $n(a) = a\bar a = c_1^2-c_2^2-c_3^2+c_4^2-c_5^2+c_6^2+c_7^2-c_8^2)$. Here $c_1,\dots,c_8 \in \mathbb C$ are the coefficients of $a \in O$.

The Jordan algebra $H_3(O,\mathbb C)$ consists of those $3 \times 3$ matrices $A = (a_{ij})$ which satisfy $a_{ij} = \bar a_{ji}$. The multiplication is given by $A \bullet B = \frac{AB+BA}{2}$. The algebra $H_3(O,\mathbb C)$ is 27-dimensional over $\mathbb C$ and 54-dimensional over $\mathbb R$. It is not hard to check that the complex form \eqref{trace_form} on $H_3(O,\mathbb C)$ is given by $t_{\mathbb C}(A) = 9(a_{11}+a_{22}+a_{33}) = 9tr\,A$. The linear form $tr\,A$ is called the {\sl generic trace} \cite[p.233]{Jacobson68}, and gives rise to a symmetric {\sl generic trace bilinear form} $t(A,B) = tr\,(A \bullet B)$ \cite[p.227]{Jacobson68}. Hence the generic trace bilinear form is related to the form \eqref{bilinear} by
\begin{equation} \label{g_t}
t(A,B) = \frac19g_{\mathbb C}(A,B).
\end{equation}
On $H_3(O,\mathbb C)$ there exists a cubic polynomial similar to the determinant, the {\sl generic norm}. It is explicitly given by \cite[eq.(50), p.232]{Jacobson68}
\[ \det A = a_{11}a_{22}a_{33} - a_{11}n(a_{23}) - a_{22}n(a_{31}) - a_{33}n(a_{12}) + C((a_{12}a_{23})a_{31}),
\]
where $C(a) = a+\bar a$ is twice the complex part of $a$. The generic trace bilinear form and the generic norm are related by \cite[eq.(68'), p.243]{Jacobson68}
\[ D_U\log \det A = t(A^{-1},U).
\]
Comparing this with \eqref{der_log_det} and taking into account that $\det I = \det P_I^{\mathbb C} = 1$, we obtain by virtue of \eqref{g_t} that $\det P_A^{\mathbb C} = (\det A)^{18}$. It follows that $\omega(A) = |\det A|^{18}$.

\subsection{Real central-simple Jordan algebras}

In this subsection we provide the classification of the real central-simple Jordan algebras and compute the $\omega$-functions of their $\omega$-domains. The complexification of a real central-simple Jordan algebra is a complex simple Jordan algebra \cite[Theorem 9, p.206]{Jacobson68} and must hence be isomorphic to one of the algebras listed in the previous subsection. As in the theory of Lie algebras, one says that the real algebra is a {\sl form} of the complex one \cite[p.70]{McCrimmon}. A complex algebra may have several non-isomorphic real forms. The invertible central elements of a real central-simple Jordan algebra have the form $z = \alpha e$ for some real $\alpha \not= 0$. As in the previous subsection, it follows by Theorem \ref{form_exist} that any non-degenerate trace form $\gamma$ on $J$ must be given by $\gamma(u,v) = g(z \bullet u,v) = \alpha g(u,v)$. As in Lemma \ref{ze_level} it follows that the 1-form $\zeta$ defined by such a trace form has a global potential which is proportional to \eqref{Phi_omega}.

{\theorem \cite[pp.207--212; p.369]{Jacobson68} Let $J$ be a finite-dimensional central-simple Jordan algebra over $\mathbb R$. Then $J$ is exactly one of the following:

\begin{itemize}
\item $\mathbb R$,
\item ${\cal J}ord_m(Q_{\mathbb R})$, a real quadratic factor for $m \geq 3$,
\item $M_m(\mathbb R)$, the algebra of real $m \times m$ matrices for $m \geq 3$,
\item $M_m(\mathbb H)$, the algebra of quaternionic $m \times m$ matrices for $m \geq 2$,
\item $S_m(\mathbb R,\Gamma)$, the twisted algebra of real symmetric $m \times m$ matrices for $m \geq 3$,
\item $H_m(\mathbb C,\Gamma)$, the twisted algebra of complex Hermitian $m \times m$ matrices for $m \geq 3$,
\item $H_m(\mathbb H,\Gamma)$, the twisted algebra of quaternionic Hermitian $m \times m$ matrices for $m \geq 3$,
\item $H_m(Q,\mathbb R)$, the algebra of $m \times m$ Hermitian matrices with entries being split quaternions over $\mathbb R$ for $m \geq 3$,
\item $SH_m(\mathbb H)$, the algebra of $m \times m$ skew-Hermitian quaternionic matrices for $m \geq 2$,
\item $H_3(\mathbb O,\Gamma)$, the twisted algebra of octonionic Hermitian $3 \times 3$ matrices,
\item $H_3(O,\mathbb R)$, the algebra of $3 \times 3$ Hermitian matrices with entries being split octonions over $\mathbb R$.
\end{itemize} }

\subsubsection{Ground field $\mathbb R$}

The Jordan multiplication in $\mathbb R$ coincides with the usual multiplication. An element is invertible if and only if it is nonzero. Hence the set of invertible elements has two connected components, namely the open half-rays. These components are mutually isomorphic. The complexification of $\mathbb R$ is $\mathbb C$. We have $L_x = x$, $P_x = x^2$, and hence $\omega(x) = |x|$.

\subsubsection{Quadratic factor ${\cal J}ord_m(Q_{\mathbb R})$}

The real quadratic factors are forms of the complex quadratic factor. They are defined in the same way as the complex quadratic factor, with the difference that $Q$ is a symmetric non-degenerate form on $\mathbb R^m$. Since $Q$ must evaluate to 1 on the unit element, it cannot be negative definite. All other signatures can occur and yield non-isomorphic Jordan algebras. Analogous to the complex case, we obtain $\det P_x = (x^TQx)^m$, and a point $x$ is invertible if and only if $x^TQx \not= 0$. The $\omega$-function is given by $\omega(x) = |x^TQx|^{m/2}$.

\subsubsection{Full real matrices $M_m(\mathbb R)$}

The Jordan algebra $M_m(\mathbb R)$ is the space of real $m \times m$ matrices equipped with the multiplication $A \bullet B = \frac{AB+BA}{2}$ \cite[p.58]{McCrimmon}. Its complexification is $M_m(\mathbb C)$. As in the complex case, $P_AB = ABA$, and $A$ is invertible in the Jordan algebra if and only if it is invertible as a matrix. Hence the set of invertible elements has two connection components, with positive and negative determinant, respectively. These components are mutually isomorphic.

Assume that the matrix $A$ is diagonalizable, with real eigenvalues $\lambda_1,\dots,\lambda_m$. By the same reasoning as in the complex case, we obtain $\det P_A = (\det A)^{2m}$. Since the set of diagonalizable matrices with real eigenvalues is open and $\det P_A$ is polynomial in $A$, this formula is valid for all matrices $A$. We then obtain $\omega(A) = |\det A|^m$.

\subsubsection{Full quaternionic matrices $M_m(\mathbb H)$}

The Jordan algebra $M_m(\mathbb H)$ is the space of quaternionic $m \times m$ matrices. A quaternion can be represented as a complex $2 \times 2$ matrix
\begin{equation} \label{quat_rep}
\begin{pmatrix} z & w \\ -\bar w & \bar z \end{pmatrix},
\end{equation}
hence $M_m(\mathbb H)$ can be represented as the subspace of complex $2m \times 2m$ matrices of the form
\begin{equation} \label{quat_form}
S = \begin{pmatrix} Z & W \\ -\bar W & \bar Z \end{pmatrix}.
\end{equation}
The Jordan multiplication on this space is given by $S \bullet T = \frac{ST+TS}{2}$ \cite[p.58]{McCrimmon}. The complexification of the algebra is isomorphic to $M_{2m}(\mathbb C)$. Again we have $P_ST = STS$, and $S$ is invertible in the Jordan algebra if and only if it is invertible as a matrix.

Note that if $S$ has an eigenvector $u = \begin{pmatrix} v \\ w \end{pmatrix}$ with eigenvalue $\lambda$, then $\tilde u = \begin{pmatrix} -\bar w \\ \bar v \end{pmatrix}$ is an eigenvector with eigenvalue $\bar\lambda$. Hence $\det S$ is real and non-negative. Moreover, since $S^T$ has the same form as $S$, this relation holds also for $S^T$. Assume that the matrix $S$ is diagonalizable, with eigenvalues $\lambda_1,\dots,\lambda_m,\bar \lambda_1,\dots,\bar \lambda_m$. Let the corresponding eigenvectors of $S$ be $u_1,\dots,u_m,\tilde u_1,\dots,\tilde u_m$ and those of $S^T$ be $y_1,\dots,y_m,\tilde y_1,\dots,\tilde y_m$. Then the matrices $T_{1,kl} = u_ky_l^T+\tilde u_k\tilde y_l^T$, $T_{2,kl} = i(u_ky_l^T - \tilde u_k\tilde y_l^T)$, $T_{3,kl} = \tilde u_ky_l^T - u_k\tilde y_l^T$, $T_{4,kl} = i(\tilde u_ky_l^T + u_k\tilde y_l^T)$, $k,l = 1,\dots,m$, constitute a basis of $M_m(\mathbb H)$ over $\mathbb R$. Note that $T_{1,kl}\pm iT_{2,kl}$ are eigenvectors of the operator $P_S$ with eigenvalues $\bar\lambda_k\bar\lambda_l,\lambda_k\lambda_l$, respectively, and $T_{3,kl}\pm iT_{4,kl}$ are eigenvectors with eigenvalues $\lambda_k\bar\lambda_l,\bar\lambda_k\lambda_l$, respectively. Thus we obtain
\[ \det P_S = \prod_{k,l=1}^m (\bar\lambda_k\bar\lambda_l)(\lambda_k\lambda_l)(\lambda_k\bar\lambda_l)(\bar\lambda_k\lambda_l) = \prod_{k,l=1}^m |\lambda_k|^4|\lambda_l|^4 = \prod_{k=1}^m |\lambda_k|^{8m} = (\det S)^{4m}.
\]
As above, this formula must be valid for all matrices $S$. We then obtain $\omega(S) = (\det S)^{2m}$.

\subsubsection{Twisted real symmetric matrices $S_m(\mathbb R,\Gamma)$}

The algebra $S_m(\mathbb R)$ of real symmetric $m \times m$ matrices with multiplication $A \bullet B = \frac{AB+BA}{2}$ is a formally real Jordan algebra. Its complexification is the algebra of complex symmetric $m \times m$ matrices $S_m(\mathbb C)$. The twisted algebras $S_m(\mathbb R,\Gamma)$ are forms of $S_m(\mathbb C)$ which are not necessarily isomorphic, but isotopic to $S_m(\mathbb R)$. Here $\Gamma$ is a diagonal matrix with diagonal elements $\pm1$, and $S_m(\mathbb R,\Gamma)$ is defined as the $\Gamma$-isotope of $S_m(\mathbb R,\Gamma)$ \cite[pp.72--73]{McCrimmon}. Thus the real vector space underlying $S_m(\mathbb R,\Gamma)$ is the same as for $S_m(\mathbb R)$, but the Jordan product is defined differently. Since $S_m(\mathbb R,\Gamma)$ is isomorphic to $S_m(\mathbb R,-\Gamma)$, one can assume that $\Gamma$ has not less positive elements than negative elements on the diagonal.

Let us first compute $\det P_A$ for the algebra $S_m(\mathbb R)$. As above, $P_AB = ABA$, and as for $S_m(\mathbb C)$ we obtain $\det P_A = (\det A)^{m+1}$. It follows that $\det P_{\Gamma} = \pm1$, depending on $m$ and the signature of $\Gamma$. By virtue of \eqref{P_prod} we then obtain that in $S_m(\mathbb R,\Gamma)$ the determinant of the operator $P_A$ is given by $\pm(\det A)^{m+1}$. Thus $A$ is invertible in $S_m(\mathbb R,\Gamma)$ if and only if it is invertible as a matrix. The connected component of the unit element $\Gamma$ in the set of invertible elements is then the set of matrices which have the same signature as $\Gamma$, and the $\omega$-function is given by $\omega(A) = |\det A|^{(m+1)/2}$.

\subsubsection{Twisted complex Hermitian matrices $H_m(\mathbb C,\Gamma)$}

The case of $H_m(\mathbb C,\Gamma)$ is similar to that of $S_m(\mathbb R,\Gamma)$. The algebras $H_m(\mathbb C,\Gamma)$ are real forms of the algebra $M_m(\mathbb C)$ of full complex $m \times m$ matrices, and they are isotopic to the formally real algebra $H_m(\mathbb C)$ of complex Hermitian $m \times m$ matrices with Jordan product $A \bullet B = \frac{AB+BA}{2}$. If $v_1,\dots,v_m \in \mathbb C^m$ are eigenvectors of $A$ with eigenvalues $\lambda_1,\dots,\lambda_m \in \mathbb R$, then $B_{ij} = v_iv_j^* + v_jv_i^*$ is an eigenvector of $P_A: B \mapsto ABA$ in $H_m(\mathbb C)$ with eigenvalue $\lambda_i\lambda_j$. Since the matrices $B_{ij}$, $i,j = 1,\dots,m$, constitute a basis of $H_m(\mathbb C)$ over $\mathbb R$, we get
\[ \det P_A = \prod_{i,j = 1}^m \lambda_i\lambda_j = \prod_{i=1}^m \lambda_i^{2m} = (\det A)^{2m}.
\]
It follows that $\det P_{\Gamma} = 1$. Repeating the above reasoning, we then get $\det P_A = (\det A)^{2m}$ in $H_m(\mathbb C,\Gamma)$, and $\omega(A) = |\det A|^m$.

\subsubsection{Twisted complex quaternionic matrices $H_m(\mathbb H,\Gamma)$}

The case of $H_m(\mathbb H,\Gamma)$ is similar to that of $H_m(\mathbb C,\Gamma)$. The algebras $H_m(\mathbb H,\Gamma)$ are real forms of the algebra $H_m(Q,\mathbb C)$ of Hermitian complex split quaternionic $m \times m$ matrices, and they are isotopic to the formally real algebra $H_m(\mathbb H)$ of quaternionic Hermitian $m \times m$ matrices with Jordan product $A \bullet B = \frac{AB+BA}{2}$. If we represent the quaternions by $2 \times 2$ complex matrices as in \eqref{quat_rep}, then $H_m(\mathbb H)$ will be represented by complex Hermitian $2m \times 2m$ matrices $S$ of the form \eqref{quat_form}.

Let us compute $\det P_S$ for the algebra $H_m(\mathbb H)$. As in the case of $M_m(\mathbb H)$, let $u_1,\dots,u_m,\tilde u_1,\dots,\tilde u_m \in \mathbb C^{2m}$ be eigenvectors of $S$ with eigenvalues $\lambda_1,\dots,\lambda_m,\lambda_1,\dots,\lambda_m$. Note that the eigenvalues are real, because the matrix $S$ is Hermitian. Then the matrices $T_k = u_ku_k^* + \tilde u_k\tilde u_k^*$, $k = 1,\dots,m$, $T_{1,kl} = u_ku_l^* + \tilde u_k\tilde u_l^* + u_lu_k^* + \tilde u_l\tilde u_k^*$, $T_{2,kl} = i(u_ku_l^* - \tilde u_k\tilde u_l^* - u_lu_k^* + \tilde u_l\tilde u_k^*)$, $T_{3,kl} = \tilde u_ku_l^* - u_k\tilde y_l^* + u_l\tilde u_k^* - \tilde u_lu_k^*$, $T_{4,kl} = i(\tilde u_ku_l^* + u_k\tilde y_l^* - u_l\tilde u_k^* - \tilde u_lu_k^*)$, $1 \leq k < l \leq m$, constitute a basis of $H_m^{\mathbb R}(\mathbb H)$. The matrices $T_k,T_{j,kl}$ are eigenvectors of $P_S: T \mapsto STS$ with eigenvalues $\lambda_k^2,\lambda_k\lambda_l$, respectively, for $j = 1,2,3,4$. It follows that
\[ \det P_S = \left( \prod_{k=1}^m \lambda_k^2 \right) \left( \prod_{k < l} \lambda_k^4\lambda_l^4 \right) = \prod_{k=1}^m \lambda_k^{2+4(m-1)} = (\det S)^{2m-1}.
\]
Since $\det S \geq 0$, we have $\det P_{\Gamma} = 1$ and hence as above $\det P_S = (\det S)^{2m-1}$ also in $H_m(\mathbb C,\Gamma)$. Thus $\omega(S) = (\det S)^{m-1/2}$.

\subsubsection{Split quaternionic matrices $H_m(Q,\mathbb R)$}

Similar to the case $H_m(Q,\mathbb C)$ we represent the split quaternions over $\mathbb R$ by real $2 \times 2$ matrices with conjugation \eqref{split_quat_conj}, and the algebra $H_m(Q,\mathbb R)$ by the algebra $A_{2m}(\mathbb R)$ of real skew-symmetric $2m \times 2m$ matrices with multiplication \eqref{split_quat_mult}. Then the operator $P_S$ acts like $T \mapsto SJTJS$. Introduce again the operator ${\bf J}$ with determinant $(-1)^{m^2}$ acting on $A_{2m}(\mathbb R)$ like $T \mapsto JTJ$. Let $v_1,\dots,v_m,\bar v_1,\dots,\bar v_m$ be the eigenvectors of an invertible skew-symmetric matrix $S$ with eigenvalues $i\lambda_1,\dots,i\lambda_m,-i\lambda_1,\dots,-i\lambda_m$, respectively, where $\lambda_k \in \mathbb R$. Then the matrices $T_k = i(v_k\bar v_k^T - \bar v_kv_k^T)$, $k = 1,\dots,m$, $T_{1,kl} = v_kv_l^T - v_lv_k^T + \bar v_k\bar v_l^T - \bar v_l\bar v_k^T$, $T_{2,kl} = i(v_kv_l^T - v_lv_k^T - \bar v_k\bar v_l^T + \bar v_l\bar v_k^T)$, $T_{3,kl} = v_k\bar v_l^T - \bar v_lv_k^T + \bar v_kv_l^T - v_l\bar v_k^T$, $T_{4,kl} = i(v_k\bar v_l^T - \bar v_lv_k^T - \bar v_kv_l^T + v_l\bar v_k^T)$, $1 \leq k < l \leq m$, constitute a basis of $A_{2m}(\mathbb R)$. Moreover, the matrices $T_k,T_{1,kl},T_{2,kl},T_{3,kl},T_{4,kl}$ are eigenvectors of $P_S \circ {\cal J}: T \mapsto STS$ with eigenvalues $-\lambda_k^2,\lambda_k\lambda_l,\lambda_k\lambda_l,-\lambda_k\lambda_l,-\lambda_k\lambda_l$, respectively. Hence
\[ \det P_S \cdot \det {\bf J} = \det (P_S \circ {\bf J}) = \prod_{k=1}^m (-\lambda_k^2) \cdot \prod_{k < l} (\lambda_k^4\lambda_l^4) = (-1)^m \prod_{k=1}^m \lambda_k^{2+4(m-1)} = (-1)^m (\det S)^{2m-1}.
\]
It follows that $\det P_S = (-1)^{m-m^2}(\det S)^{2m-1} = (\pf S)^{2(2m-1)}$ and $\omega(S) = |\pf S|^{2m-1}$.

\subsubsection{Skew-Hermitian quaternionic matrices $SH_m(\mathbb H)$}

Let $Q$ be a non-degenerate skew-Hermitian quaternionic $m \times m$ matrix. Then the multiplication $A \bullet B = \frac{AQB + BQA}{2}$ on the space of skew-Hermitian quaternionic matrices defines a Jordan algebra with unit element $Q^{-1}$. This algebra is a real form of the algebra $S_{2m}(\mathbb C)$ of complex symmetric $2m \times 2m$ matrices. The isomorphism class of the algebra does not depend on the choice of $Q$ (cf.\ \cite[Ex.5, p.211]{Jacobson68}). We shall choose $Q = iI$. If we represent the quaternions by complex $2 \times 2$ matrices \eqref{quat_rep}, then the elements of $SH_m(\mathbb H)$ can be represented by complex skew-Hermitian $2m \times 2m$ matrices $S$ of the form \eqref{quat_form}. If we define the matrix $\Lambda = \diag(iI,-iI)$, then the multiplication in $SH_m(\mathbb H)$ is given by $S \bullet T = \frac{S\Lambda T+T\Lambda S}{2}$. It follows that the operator $P_S$ acts like $T \mapsto S\Lambda T \Lambda S$.

Introduce the operator ${\bf \Lambda}: T \mapsto \Lambda T\Lambda$. It is not hard to see that $\det{\bf \Lambda} = (-1)^{m^2}$. The composition $P_S \circ {\bf \Lambda}$ acts like $T \mapsto S\Lambda^2 T \Lambda^2 S = STS$. As in the case of $M_m(\mathbb H)$, let $u_1,\dots,u_m,\tilde u_1,\dots,\tilde u_m \in \mathbb C^{2m}$ be eigenvectors of $S$ with eigenvalues $i\lambda_1,\dots,i\lambda_m,-i\lambda_1,\dots,-i\lambda_m$, where $\lambda_k \in \mathbb R$ because $S$ is skew-Hermitian. Then the matrices $B_{1,k} = i(u_ku_k^* - \tilde u_k\tilde u_k^*)$, $B_{2,k} = \tilde u_ku_k^* - u_k\tilde u_k^*$, $B_{3,k} = i(\tilde u_ku_k^* + u_k\tilde u_k^*)$, $k = 1,\dots,m$, $B_{1,kl} = u_ku_l^* + \tilde u_k\tilde u_l^* - u_lu_k^* - \tilde u_l\tilde u_k^*$, $B_{2,kl} = i(u_ku_l^* - \tilde u_k\tilde u_l^* + u_lu_k^* - \tilde u_l\tilde u_k^*)$, $B_{3,kl} = \tilde u_ku_l^* - u_k\tilde u_l^* - u_l\tilde u_k^* + \tilde u_lu_k^*$, $B_{4,kl} = i(\tilde u_ku_l^* + u_k\tilde y_l^* + u_l\tilde u_k^* + \tilde u_lu_k^*)$, $1 \leq k < l \leq m$, constitute a basis of $SH_m(\mathbb H)$. These matrices are eigenvectors of $P_S \circ {\bf \Lambda}$ with eigenvalues $-\lambda_k^2,\lambda_k^2,\lambda_k^2,-\lambda_k\lambda_l,-\lambda_k\lambda_l,\lambda_k\lambda_l,\lambda_k\lambda_l$, respectively. It follows that
\[ \det P_S \cdot \det {\bf \Lambda} = \det (P_S \circ {\bf \Lambda}) = \prod_{k=1}^m (-\lambda_k^6) \cdot \prod_{k < l} (\lambda_k^4\lambda_l^4) = (-1)^m \prod_{k=1}^m \lambda_k^{6+4(m-1)} = (-1)^m (\det S)^{2m+1}.
\]
Hence $\det P_S = (-1)^{m-m^2}(\det S)^{2m+1} = (\det S)^{2m+1}$ and $\omega(S) = (\det S)^{m+1/2}$.

\subsubsection{Twisted octonionic Hermitian matrices $H_3(\mathbb O,\Gamma)$}

The case of $H_3(\mathbb O,\Gamma)$ is similar to that of the other twisted Hermitian matrix algebras. There are two non-isomorphic algebras $H_3(\mathbb O,\Gamma)$, for $\Gamma = I$ and for $\Gamma = \diag(1,1,-1)$. For $\Gamma = I$ we obtain the formally real algebra $H_3(\mathbb O)$ of $3 \times 3$ octonionic Hermitian matrices with the usual Jordan product $A \bullet B = \frac{AB+BA}{2}$. Both algebras $H_3(\mathbb O,\Gamma)$ are real forms of the split octonion algebra $H_3(O,\mathbb C)$, and are mutually isotopic.

By the same arguments as in Subsection \ref{subsec_split_oct} we get $\det P_A = (\det A)^{18}$, where
\begin{equation} \label{det_octonion}
\det A = a_{11}a_{22}a_{33} - a_{11}n(a_{23}) - a_{22}n(a_{31}) - a_{33}n(a_{12}) + 2Re((a_{12}a_{23})a_{31}).
\end{equation}
Here $n(a) = a\bar a = |a|^2$. It follows that $\omega(A) = |\det A|^9$.

\subsubsection{Split octonionic matrices $H_3(O,\mathbb R)$}

The real split octonions are defined in the same manner as the complex split octonions in Subsection \ref{subsec_split_oct}, except that the coefficients $c_1,\dots,c_8$ are required to be real. The conjugate $\bar a$ and the norm $n(a)$ of an element $a$ are also defined by the same formulas as in the complex case.

The Jordan algebra $H_3(O,\mathbb R)$ of $3 \times 3$ Hermitian matrices with real split octonionic entries and multiplication $A \bullet B = \frac{AB+BA}{2}$ is 27-dimensional over $\mathbb R$, with $H_3(O,\mathbb C)$ being its complexification. By the same arguments as in Subsection \ref{subsec_split_oct} we get $\det P_A = (\det A)^{18}$, where $\det A$ is given by \eqref{det_octonion}. It again follows that $\omega(A) = |\det A|^9$.

\subsection{Final classification}

In this subsection we provide a complete classification of proper affine hyperspheres in $\mathbb R^n$ with parallel cubic form. Actually, we classify all centro-affine hypersurface immersions with parallel cubic form which are associated to a semi-simple Jordan algebra $J$. By Lemma \ref{sum_zeta_Phi} the 1-form $\zeta$ defining such an immersion locally has a potential which can be represented as a sum of potentials $\Phi$ defined on the individual simple factors of $J$. In the last two subsections we computed these potentials $\Phi$ for all simple real Jordan algebras. Note that we are not really interested in the Jordan algebra $J$, but only in the underlying vector space and in the functions $\Phi$ and $\omega$ on it.

Let us provide these data in the form of a table. In the first column we list the vector space, and in the second column its dimension. Here $M_m,S_m,A_m,H_m,SH_m$ stands for full, symmetric, skew-symmetric, Hermitian and skew-Hermitian matrices of size $m \times m$, respectively. Most of the classes of real simple Jordan algebras constitute infinite series parameterized by an integer. We give the range of this parameter in the third column. In the fourth column we give an expression for the local potential $\Phi$, parameterized by a nonzero complex number $c$ for complex Jordan algebras, and by a nonzero real number $\alpha$ for real central-simple algebras. In the last two columns we provide the $\omega$-function of the corresponding $\omega$-domains and a description of the affine hyperspheres associated with these domains. The constants in the last column are assumed to be nonzero. Note that in the case of a matrix space over the quaternions $\mathbb H$, the matrix $S$ is the complex representation \eqref{quat_form} of the quaternionic matrix and has twice the size. In the row corresponding to the vector space $\mathbb R^m$, $Q$ denotes a non-degenerate quadratic form on $\mathbb R^m$.

\medskip

\begin{tabular}{|c|c|c|c|c|c|}
\hline
vector space & real dimension & range & $\Phi$ & $\omega$ & affine sphere \\
\hline
\hline
$\mathbb C$ & 2 & & $Re(c\log x)$ & $|x|^2$ & $|x| = const$ \\
\hline
$\mathbb C^m$ & $2m$ & $m \geq 3$ & $Re(c\log x^Tx)$ & $|x^Tx|^m$ & $|x^Tx| = const$ \\
\hline
$S_m(\mathbb C)$ & $m(m+1)$ & $m \geq 3$ & $Re(c\log\det A)$ & $|\det A|^{m+1}$ & $|\det A| = const$ \\
\hline
$M_m(\mathbb C)$ & $2m^2$ & $m \geq 3$ & $Re(c\log\det A)$ & $|\det A|^{2m}$ & $|\det A| = const$ \\
\hline
$A_{2m}(\mathbb C)$ & $2m(2m-1)$ & $m \geq 3$ & $Re(c\log\pf A)$ & $|\pf A|^{2(2m-1)}$ & $|\pf A| = const$ \\
\hline
$H_3(O,\mathbb C)$ & 54 & & $Re(c\log\det A)$ & $|\det A|^{18}$ & $|\det A| = const$ \\
\hline
\hline
$\mathbb R$ & 1 & & $\alpha\log|x|$ & $|x|$ & point \\
\hline
$\mathbb R^m$ & $m$ & $m \geq 3$ & $\alpha\log|x^TQx|$ & $|x^TQx|^{m/2}$ & quadric \\
\hline
$M_m(\mathbb R)$ & $m^2$ & $m \geq 3$ & $\alpha\log|\det A|$ & $|\det A|^m$ & $\det A = const$ \\
\hline
$M_m(\mathbb H)$ & $4m^2$ & $m \geq 2$ & $\alpha\log\det S$ & $(\det S)^{2m}$ & $\det S = const$ \\
\hline
$S_m(\mathbb R)$ & $\frac{m(m+1)}{2}$ & $m \geq 3$ & $\alpha\log|\det A|$ & $|\det A|^{(m+1)/2}$ & $\det A = const$ \\
\hline
$H_m(\mathbb C)$ & $m^2$ & $m \geq 3$ & $\alpha\log|\det A|$ & $|\det A|^m$ & $\det A = const$ \\
\hline
$H_m(\mathbb H)$ & $m(2m-1)$ & $m \geq 3$ & $\alpha\log\det S$ & $(\det S)^{m-1/2}$ & $\det S = const$ \\
\hline
$A_{2m}(\mathbb R)$ & $m(2m-1)$ & $m \geq 3$ & $\alpha\log|\pf A|$ & $|\pf A|^{2m-1}$ & $\pf A = const$ \\
\hline
$SH_m(\mathbb H)$ & $m(2m+1)$ & $m \geq 2$ & $\alpha\log\det S$ & $(\det S)^{m+1/2}$ & $\det S = const$ \\
\hline
$H_3(\mathbb O)$ & 27 & & $\alpha\log|\det A|$ & $|\det A|^9$ & $\det A = const$ \\
\hline
$H_3(O,\mathbb R)$ & 27 & & $\alpha\log|\det A|$ & $|\det A|^9$ & $\det A = const$ \\
\hline
\end{tabular}

\medskip

{\theorem \label{main_centro} Let $M \subset \mathbb R^n$ be a connected non-degenerate centro-affine hypersurface with parallel cubic form, such that the real unital Jordan algebra defined by $M$ as in Theorem \ref{main1} is semi-simple.

Then $\mathbb R^n$ can be decomposed into a direct product of vector spaces $V_1,\dots,V_r$, where each of the $V_k$ is a space indicated in the first column of the above table. There locally exists a scalar function $\Phi = \sum_{k=1}^r \Phi_k$ on $\mathbb R^n$ such that $M$ can be described as a level surface of $\Phi$, $\Phi$ is logarithmically homogeneous with parameter $\nu = 1$, and each of the $\Phi_k$ is locally defined on $V_k$ and has the form indicated in the fourth column of the table for some nonzero $c \in \mathbb C$ or $\alpha \in \mathbb R$, respectively.

On the other hand, the level surfaces of any such sum $\Phi = \sum_{k=1}^r \Phi_k$ are centro-affine hypersurfaces with parallel cubic form, provided the logarithmic homogeneity constant of $\Phi$ is nonzero. }

\begin{proof}
The first part follows from Theorem \ref{th_main3}, Lemma \ref{sum_zeta_Phi}, and the completeness of the above classification of real simple Jordan algebras. The second part of the theorem follows from Lemma \ref{reverse_ca}.
\end{proof}

{\theorem \label{main_sphere} Let $M \subset \mathbb R^n$ be a proper affine hypersphere with parallel cubic form and with center in the origin. Then $\mathbb R^n$ can be decomposed into a direct product of vector spaces $V_1,\dots,V_r$, where each of the $V_k$ is a space indicated in the first column of the above table, and $M$ is a Calabi product of proper affine hyperspheres $M_k \subset V_k$ which have the form indicated in the last column of the above table.

On the other hand, all affine hyperspheres listed in the last column of the table have parallel cubic form. }

\begin{proof}
The first part follows from Theorem \ref{th_AHS2}, Lemma \ref{lem_Calabi}, and the completeness of the above classification. The second part follows from Theorem \ref{S1}.
\end{proof}

\subsection{Examples}

In this subsection we give some concrete examples of centro-affine hypersurface immersions with parallel cubic form.

\smallskip

Consider the hypersurfaces defined by the Jordan algebra $\mathbb C$ and the complex number $c = 1+i\beta$, $\beta \not= 0$. These hypersurfaces are locally given by the relation $Re(c\log x) = const$. Let $x = r\exp(i\varphi)$. Then we obtain the relation $Re((1+i\beta)(\log r+i\varphi)) = \log r - \beta\varphi = const$. It follows that the hypersurfaces are logarithmic spirals. Since the origin is contained in the closure of such a spiral, it is affine complete, but not Euclidean complete.

\smallskip

Consider the hypersurfaces defined by the product $\mathbb C \times \mathbb R$ and the non-zero numbers $c = i\beta$, $\alpha$. Let $(z = r\exp(i\varphi),x)$ be the coordinates in $\mathbb C \times \mathbb R$. Then the hypersurfaces are given by $Re(c\log z) + \alpha\log|x| = -\beta\varphi + \alpha\log|x| = const$. In Cartesian coordinates $x_1,x_2,x_3$ these surfaces are locally given by $x_3 = const \cdot \exp(\frac{\beta}{\alpha}\arctan\frac{x_2}{x_1})$. Although the Jordan algebra defined by these hypersurfaces is decomposable, the surfaces themselves cannot be decomposed into a product of lower-dimensional centro-affine hypersurfaces with parallel cubic form.

\smallskip

Consider the hypersurfaces defined by the product $\mathbb C \times \mathbb C$ and the non-zero numbers $c_k = a_k+ib_k$, $k = 1,2$, where $a_1+a_2 \not= 0$. Let $(z_1 = r_1\exp(i\varphi_1),z_2 = r_2\exp(i\varphi_2))$ be the coordinates in $\mathbb C \times \mathbb C$. Then the hypersurfaces are locally given by $Re(c_1\log z_1 + c_2\log z_2) = a_1\log r_1 + a_2\log r_2 - b_1\varphi_1 - b_2\varphi_2 = const$. For $\rho_1,\rho_2 > 0$, consider the torus $T_{\rho_1,\rho_2} \subset \mathbb C \times \mathbb C$ given by the equations $r_k = \rho_k$, $k = 1,2$. The intersection of each of the hypersurfaces with $T_{\rho_1,\rho_2}$ is given by $b_1\varphi_1 + b_2\varphi_2 = const$. If now $b_1,b_2$ have an irrational ratio, then each of the hypersurfaces has an intersection with $T_{\rho_1,\rho_2}$ which is dense in $T_{\rho_1,\rho_2}$. In this case the hypersurface itself is dense in $\mathbb C \times \mathbb C$ and cannot be an embedded submanifold.

\subsubsection{Algebras which are not semi-simple} \label{subs_non_semi_simple}

In this subsection we give a counter-example to the claim in \cite[item (v), pp.71--72]{Koecher99} and describe a family of real unital Jordan algebras with non-degenerate trace form which are not semi-simple.

Consider the $n$-dimensional vector space $J$ of real univariate polynomials $p(t)$ of degree not greater than $n-1$. On this space we can define a commutative, associative multiplication $\bullet$ by setting $p \bullet q = (p \cdot q) \mod t^n$. This multiplication turns $J$ into a real Jordan algebra with unit element $e$ given by $e(t) \equiv 1$. The inverse $p^{-1}$ of a polynomial $p \in J$ is given by the truncation of the Taylor series of the inverse $\frac{1}{p(t)}$ around $t = 0$ and exists if and only if $p(0) \not= 0$. From \cite[Theorem III.6, p.60]{Koecher99} it follows that the algebra $J$ is not semi-simple, since every polynomial $p \in J$ satisfying $p(0) = 0$ is nilpotent and the subspace of these polynomials is an ideal.

Consider the symmetric bilinear form $\gamma$ on $J$ given by $\gamma(p,q) = -(p \bullet q)(1)$. Since $J$ is associative, we have $(p \bullet q) \bullet r = p \bullet (q \bullet r)$ and hence also $\gamma(p \bullet q,r) = \gamma(p,q \bullet r)$ for all $p,q,r \in J$. Hence $\gamma$ is a trace form. Let $p = \sum_{k=0}^{n-1} p_k t^k \in J$ be a polynomial such that $\gamma(p,q) = 0$ for all $q \in J$. Inserting $q = t^k$ for $k = n-1,n-2,\dots,0$, we obtain $p_l = 0$ for $l = 0,1,\dots,n-1$, and $p = 0$. Hence $\gamma$ is non-degenerate.

Let $\zeta$ be the form \eqref{zeta_def} defined by $\gamma$. At an invertible polynomial $p \in J$ we have $\zeta(u) = -\gamma(u,p^{-1}) = (u \bullet p^{-1})(1)$. For invertible $p \in J$, define the polynomial $\log p \in J$ by the truncation of the Taylor series of $\log |p(t)|$ around $t = 0$. This Taylor series is given by
\begin{eqnarray} \label{log_Taylor}
\log|p_0+p_1t+p_2t^2+\dots| &=& \log|p_0| + \frac{p_1}{|p_0|}t + \left( - \frac12\frac{p_1^2}{p_0^2} + \frac{p_2}{|p_0|} \right)t^2 + \left( \frac13\frac{p_1^3}{|p_0|^3} - \frac{p_1p_2}{p_0^2} + \frac{p_3}{|p_0|} \right)t^3 \nonumber \\ && + \left(- \frac14\frac{p_1^4}{p_0^4} + \frac{p_1^2p_2}{|p_0|^3} - \frac{p_1p_3}{p_0^2} - \frac12\frac{p_2^2}{p_0^2} + \frac{p_4}{|p_0|} \right)t^4 + \dots
\end{eqnarray}
Since the derivations of $\log |p(t)|$ with respect to $t$ and with respect to the coefficients $p_k$ of $p$ are interchangeable, we have that $D_u((\log p)(t)) = (u \bullet p^{-1})(t)$ for all $t$. In particular, for $t = 1$ we obtain $D_u((\log p)(1)) = (u \bullet p^{-1})(1) = \zeta(u)$. Thus $F(p) = (\log p)(1)$ is a potential of the form $\zeta$, and the maximal integral manifolds of the kernel $\Delta$ of $\zeta$ are given by the relation $(\log p)(1) = const$.

The integral hypersurface through $e$ is then given by $(\log p)(1) = 0$. From \eqref{log_Taylor} we obtain the following explicit expressions of this hypersurface for $n = 2,3,4,5$:
\begin{eqnarray*}
p_1 &=& -p_0\log p_0, \\
p_2 &=& -p_0\log p_0 - p_1 + \frac12\frac{p_1^2}{p_0}, \\
p_3 &=& -p_0\log p_0 - p_1 + \frac12\frac{p_1^2}{p_0} - p_2 - \frac13\frac{p_1^3}{p_0^2} + \frac{p_1p_2}{p_0}, \\
p_4 &=& -p_0\log p_0 - p_1 + \frac12\frac{p_1^2}{p_0} - p_2 - \frac13\frac{p_1^3}{p_0^2} + \frac{p_1p_2}{p_0} - p_3 + \frac14\frac{p_1^4}{p_0^3} - \frac{p_1^2p_2}{p_0^2} + \frac{p_1p_3}{p_0} + \frac12\frac{p_2^2}{p_0}.
\end{eqnarray*}

The same construction can be carried out with spaces of multivariate polynomials, yielding a plethora of examples of centro-affine hypersurfaces with parallel cubic form which are associated to a Jordan algebra which is not semi-simple.

\section{Other classes of immersions} \label{sec_other}

In this section we show that our algebraic method is applicable also to other notions of affine hypersurface immersions with parallel cubic form or parallel difference tensor. We do not attempt a classification for these cases. We indicate, however, which classes of algebras arise from the different classes of hypersurface immersions with parallel cubic form. We consider the cases mentioned in the introduction and several other cases. Since improper affine hyperspheres arise in all of the cases mentioned in the introduction, we shall first consider graph immersions and improper affine hyperspheres in general.

\subsection{Improper affine hyperspheres}

By \cite[p.47]{NomizuSasaki} an improper affine hypersphere $M \subset \mathbb R^{n+1}$ can be represented as a graph immersion with graph function $F: \mathbb R^n \to \mathbb R$ satisfying $\det F'' \equiv \pm1$. Then the Christoffel symbols of the affine connection $\nabla$ on $M$ vanish and the affine metric $h$ is given by the Hessian metric $F''$. The Christoffel symbols of the Levi-Civita connection $\hat\nabla$ of $h$ are then given by \eqref{Christoffel} and the difference tensor $K = \nabla - \hat\nabla$ by \eqref{K_def}. At $y \in M$ the tensor $K$ defines a commutative multiplication $\bullet$ on $T_yM$ by $u \bullet v = K(u,v)$. Equipped with this multiplication, $T_yM$ becomes a commutative real algebra $A$. Moreover, by the symmetry of the third derivative $F'''$ the symmetric bilinear form $\gamma$ defined by the Hessian $F''$ on $T_yM$ satisfies condition \eqref{gamma_sym} and is hence a non-degenerate trace form. While this holds for general non-degenerate graph immersions, the property of being an affine hypersphere imposes additional conditions. Vanishing of the trace of $K$ implies that the operator $L_u$ of multiplication with the element $u$ has zero trace for all $u \in A$. In particular, $A$ cannot have a unit element.

We shall now consider different cases in more detail.

\subsubsection{Blaschke immersions with $\nabla C = 0$}

Consider non-degenerate Blaschke immersions whose cubic form is parallel with respect to the affine connection. In this case the immersion is either a quadric or a graph immersion with a cubic polynomial as graph function \cite{Vrancken88}. We shall consider the latter case. Differentiating the relation $\det F'' = const$ twice at $y \in M$ and using that the fourth derivatives of $F$ vanish identically, we obtain the condition $tr\,(L_uL_v) = 0$ for all $u,v \in A$ \cite[Theorem 2, 4)]{BNS90}. This condition is equivalent to the condition $tr\,L_u^2 = 0$ for all $u \in A$, from which it can be obtained by polarization. In general, the $k$-th derivative of the relation $\det F'' = const$ is equivalent to the condition $tr\,L_u^k = 0$ for all $u \in A$. Thus, under the assumption that $F$ is a cubic polynomial, the relation $\det F'' = const$ is equivalent to the condition that $L_u$ is nilpotent for all $u \in A$. In particular, $A$ is a nilalgebra, i.e., it consists of nilpotent elements.

On the other hand, every commutative algebra $A$ with non-degenerate trace form $\gamma$ having determinant $\pm 1$ and such that $L_u$ is nilpotent for all $u \in A$ defines an improper affine hypersphere satisfying $\nabla C = 0$ by the graph of the function $F(x) = \frac12\gamma(x,x) - \frac13\gamma(x,x^2)$.

\subsubsection{Blaschke immersions with $\nabla K = 0$}

Let us consider non-degenerate Blaschke immersions whose difference tensor is parallel with respect to the affine connection $\nabla$. These immersions are also either quadrics or improper affine hyperspheres. In the latter case the graph function is given by a polynomial, the affine metric is flat, there exists an integer $m$ such that $K_X^m = 0$ for all $X$, and $[K_X,K_Y] = 0$ for all vector fields $X,Y$ \cite{DillenVrancken98}. These conditions mean that the operator $L_u$ in the algebra $A$ is nilpotent for all $u \in A$ and that $L_u,L_v$ commute for all $u,v \in A$. We then have
\[ x \bullet (y \bullet z) = x \bullet (z \bullet y) = z \bullet (x \bullet y) = (x \bullet y) \bullet z,
\]
where the first and the last equality follow from commutativity, and the second one because $[L_x,L_z] = 0$. Hence $A$ is associative. But then $A$ is nilpotent, i.e., there exists an integer $m'$ such that the product of $m'$ arbitrary elements is zero.

On the other hand, every commutative, associative, nilpotent algebra $A$ with non-degenerate trace form $\gamma$ having determinant $\pm 1$ defines an improper affine hypersphere satisfying $\nabla K = 0$ by the graph of the function $F(x) = \sum_{k=2}^{\infty} \frac{(-2)^{k-2}}{k!}\gamma(x,x^{k-1})$. Note that the sum is finite because $x$ is nilpotent.

\subsubsection{Improper hyperspheres with $\hat\nabla C = 0$}

We now consider non-degenerate Blaschke immersions whose cubic form is parallel with respect to the Levi-Civita connection $\hat\nabla$ of the affine metric. In \cite{BNS90} it was shown that such an immersion is an affine hypersphere. This hypersphere is either proper or improper. Since every proper affine hypersphere with center in the origin has affine metric proportional to the centro-affine metric, the proper hyperspheres satisfying $\hat\nabla C = 0$ are covered by Theorem \ref{main_sphere}. Let us consider the case of improper affine hyperspheres.

Note that the cubic form $C$ is given by the third derivative $F'''$ and the affine metric by the second derivative $F''$. Then \eqref{parallel_F3} implies that the condition $\hat\nabla C = 0$ is equivalent to \eqref{quasi_lin_PDE}. The integrability condition of this PDE again leads to \eqref{int_cond}, and the algebra $A$ satisfies the Jordan identity in Definition \ref{Jordan_def}. Hence $A$ is a Jordan algebra. Now, however, the form \eqref{trace_form} vanishes identically, and thus also the bilinear form \eqref{bilinear} vanishes identically. From \cite[Section III.2]{Koecher99} it follows that $A$ is a nilalgebra, and hence nilpotent \cite[pp.195--196]{Jacobson68}.

In a subsequent publication we will show that the converse is also true. Every nilpotent Jordan algebra with non-degenerate trace form $\gamma$ having determinant $\pm 1$ defines an improper affine hypersphere satisfying $\hat\nabla C = 0$ by the graph of the function $F(x) = \sum_{k=2}^{\infty} \frac{(-1)^k}{k} \gamma(x,x^{k-1})$. Again the sum is finite because $x$ is nilpotent.

\subsection{Centro-affine immersions}

We shall now consider several additional cases of centro-affine immersions. As in the main part of the paper, we define the algebra $A$ by the difference tensor \eqref{K_def}. By \eqref{pos_identity} the position vector $e$ is the unit element of the algebra $A$. Moreover, the symmetric bilinear form $\gamma$ defined by the Hessian $F''$ satisfies \eqref{gamma_symmetric} and is hence a non-degenerate trace form. By the fifth relation in Lemma \ref{FderX} we also have $\gamma(e,e) = -1$.

\subsubsection{Centro-affine immersions with $\nabla C = 0$}

Consider centro-affine immersions whose cubic form is parallel with respect to the affine connection $\nabla$. A somewhat tedious calculus shows that the condition $\nabla C = 0$ is equivalent to the quasi-linear PDE
\begin{eqnarray*}
\lefteqn{F_{,ijkl} + 3(F_{,ijl}F_{,k} + F_{,ikl}F_{,j} + F_{,jkl}F_{,i} + F_{,ijk}F_{,l}) + 2(F_{,ij}F_{,kl} + F_{,ik}F_{,jl} + F_{,jk}F_{,il}) } && \\ &+& 8(F_{,il}F_{,j}F_{,k} + F_{,jl}F_{,i}F_{,k} + F_{,kl}F_{,i}F_{,j} + F_{,jk}F_{,i}F_{,l} + F_{,ik}F_{,j}F_{,l} + F_{,ij}F_{,k}F_{,l}) + 24F_{,i}F_{,j}F_{,k}F_{,l} = 0.
\end{eqnarray*}
Differentiating with respect to $x^m$, replacing the appearing fourth derivatives, and anti-commutating the indices $l,m$ yields
\begin{eqnarray*}
\lefteqn{-(F_{,jkm}F_{,i}F_{,l} + F_{,ikm}F_{,j}F_{,l} + F_{,ijm}F_{,k}F_{,l}) + 2(F_{,jk}F_{,im}F_{,l} + F_{,ik}F_{,jm}F_{,l} + F_{,ij}F_{,km}F_{,l}) } && \\ &+& (F_{,ijl}F_{,km} + F_{,ikl}F_{,jm} + F_{,jkl}F_{,im}) = -(F_{,jkl}F_{,i}F_{,m} + F_{,ikl}F_{,j}F_{,m} + F_{,ijl}F_{,k}F_{,m}) \\ &+& 2(F_{,jk}F_{,il}F_{,m} + F_{,ik}F_{,jl}F_{,m} + F_{,ij}F_{,kl}F_{,m}) + (F_{,ijm}F_{,kl} + F_{,ikm}F_{,jl} + F_{,jkm}F_{,il}).
\end{eqnarray*}
Raising the indices $l,m$ with the pseudo-metric $F''$ and contracting $i,j,k$ with a vector $u$ yields
\[ K_{jk}^mF_{,i}e^lu^iu^ju^k + F_{,jk}u^ju^ku^me^l + K_{ij}^lu^iu^ju^m = K_{jk}^lF_{,i}e^mu^iu^ju^k + F_{,jk}u^ju^ku^le^m + K_{ij}^mu^iu^ju^l.
\]
With $v = u^2$ and by virtue of the second relation in Lemma \ref{FderX} this can be written as
\begin{equation} \label{skew_euv}
-\gamma(u,e)e^lv^m + \gamma(u,u)e^lu^m + v^lu^m = -\gamma(u,e)v^le^m + \gamma(u,u)u^le^m + u^lv^m.
\end{equation}
It follows that the vectors $e,u,u^2$ cannot be linearly independent. If $u$ is a multiple of the unit element $e$, then $u^2$ is too. Hence in any case $u^2$ is a linear combination of $u$ and $e$. It follows that $L_{u^2}$ is a linear combination of $L_u$ and the identity matrix. In particular, it commutes with $L_u$ and $A$ is a Jordan algebra. Linear dependence of $e,u,u^2$ implies that $A$ is of degree 2 and hence a quadratic factor. Using that $\gamma(u,u) = \gamma(v,e)$ it is not hard to see that this is actually sufficient for \eqref{skew_euv} to hold. In order for a quadratic factor to have a non-degenerate trace form it must either be of dimension 2, or non-degenerate and hence central-simple. In the second case the trace form $\gamma$ is proportional to the form \eqref{bilinear}.

\subsubsection{Centro-affine immersions with $\nabla K = 0$}

Consider centro-affine immersions whose difference tensor is parallel with respect to the affine connection $\nabla$. It is not hard to prove that this condition is equivalent to the relation $\nabla_lC_{ijk} = C_{ijr}h^{rs}C_{kls}$, where $h$ is the affine metric and $C$ the cubic form. A lengthy calculus shows that this is equivalent to the quasi-linear PDE
\begin{eqnarray}
\lefteqn{F_{,ijkl} - F_{,ijr}F^{,rs}F_{,kls} + 2(F_{,ik}F_{,jl} + F_{,jk}F_{,il} - F_{,ij}F_{,kl}) + F_{,ijl}F_{,k} + F_{,ikl}F_{,j} } && \label{PDE2} \\
&+& F_{,jkl}F_{,i} + F_{,ijk}F_{,l} + 4(F_{,il}F_{,j}F_{,k} + F_{,jl}F_{,i}F_{,k} + F_{,jk}F_{,i}F_{,l} + F_{,ik}F_{,j}F_{,l}) + 8F_{,i}F_{,j}F_{,k}F_{,l} = 0. \nonumber
\end{eqnarray}
Anti-symmetrizing with respect to the indices $i,k$ yields
\[ F_{,jkr}F^{,rs}F_{,ils} - F_{,ijr}F^{,rs}F_{,kls} + 4(F_{,jk}F_{,il} - F_{,ij}F_{,kl} + F_{,il}F_{,j}F_{,k} + F_{,jk}F_{,i}F_{,l} - F_{,kl}F_{,j}F_{,i} - F_{,ij}F_{,k}F_{,l}) = 0.
\]
Raising the index $j$ with the pseudo-metric $F''$ leads to
\[ K_{kr}^jK_{il}^r - K_{ir}^jK_{kl}^r + \delta_k^jF_{,il} - \delta_i^jF_{,kl} - F_{,il}e^jF_{,k} + \delta_k^jF_{,i}F_{,l} + F_{,kl}e^jF_{,i} - \delta_i^jF_{,k}F_{,l} = 0.
\]
Denote the 1-form defined by the first derivative $F'$ on $A$ by $\rho$. Contracting with vectors $u^i,v^k,w^l$ gives
\begin{equation} \label{uv_comm}
[L_v,L_u]w + (\gamma(u,w) + \rho(u)\rho(w))v - (\gamma(v,w) + \rho(v)\rho(w))u + (\gamma(v,w)\rho(u) - \gamma(u,w)\rho(v))e = 0.
\end{equation}
Symmetrizing the PDE \eqref{PDE2} leads to
\begin{eqnarray*}
\lefteqn{F_{,ijkl} - \frac13(F_{,ijr}F^{,rs}F_{,kls} + F_{,ikr}F^{,rs}F_{,jls} + F_{,ilr}F^{,rs}F_{,jks}) + \frac23(F_{,ik}F_{,jl} + F_{,jk}F_{,il} + F_{,ij}F_{,kl}) } && \\ &+& F_{,ijl}F_{,k} + F_{,ikl}F_{,j} + F_{,jkl}F_{,i} + F_{,ijk}F_{,l} + \frac83\left(F_{,il}F_{,j}F_{,k} + F_{,jl}F_{,i}F_{,k} + F_{,jk}F_{,i}F_{,l} + F_{,ik}F_{,j}F_{,l} \right. \\ &+& \left. F_{,ij}F_{,k}F_{,l} + F_{,kl}F_{,i}F_{,j}\right) + 8F_{,i}F_{,j}F_{,k}F_{,l} = 0.
\end{eqnarray*}
Differentiating this with respect to $x^m$, anti-symmetrizing the indices $l,m$, raising the index $m$ by means of the pseudo-metric $F''$, contracting the indices $i,j,k$ with a vector $u$, and the index $l$ with a vector $v$ gives after a lengthy calculation
\[ 8[L_u,L_{u^2}]v + (\rho(u)\rho(v) + \gamma(u,v))u^2 + (\gamma(u,u)\rho(v) - \gamma(v,u^2))u + (\gamma(u,u)\gamma(u,v) + \gamma(u^2,v)\rho(u))e = 0.
\]
Expressing the commutator by \eqref{uv_comm} and using that $\rho(u^2) = -\gamma(u,u)$ gives
\begin{equation} \label{second_condition}
(\gamma(u,v) + \rho(u)\rho(v))u^2 + (\gamma(u,u)\rho(v) - \gamma(u^2,v))u + (\gamma(u,u)\gamma(u,v) + \gamma(u^2,v)\rho(u))e = 0.
\end{equation}
Now by the first item in Lemma \ref{correspondence} we have $\gamma(u,v) + \rho(u)\rho(v) = 0$ for all $v$ if and only if $u$ is proportional to $e$, and hence $u^2$ is again a linear combination of the unit element $e$ and the vector $u$. As above, it follows that $A$ is a quadratic factor, and either 2-dimensional or non-degenerate. As in the previous case, one can show that this is sufficient for \eqref{uv_comm},\eqref{second_condition} to hold.

\section*{Acknowledgements}

The author would like to thank Professors Zejun Hu, An-Min Li, and Udo Simon for enlightening discussions on the subject.

\bibliography{affine_geometry,jordan,geometry,algebra}
\bibliographystyle{plain}

\end{document}